\documentclass[draft]{amsart}
\usepackage{amsmath,amssymb,amscd}

\newtheorem{theorem*}{Theorem}
\newtheorem{theorem}{Theorem}[section]
\newtheorem{proposition}[theorem]{Proposition}
\newtheorem{corollary}[theorem]{Corollary}
\newtheorem{lemma}[theorem]{Lemma}

\theoremstyle{definition}
\newtheorem{definition}[theorem]{Definition}
\newtheorem{example}[theorem]{Example}

\theoremstyle{remark}
\newtheorem*{remark}{Remark}

\numberwithin{equation}{section}

\renewcommand{\L}{\Lambda}

\renewcommand{\l}{\rho}

\newcommand{\f}{\varphi}

\newcommand{\F}{Q}

\newcommand{\mb}[1]{{\textbf {\textit#1}}}
\renewcommand{\u}[1]{ u_{{}_{\!#1}} }
\renewcommand{\v}[1]{ v_{{}_{\!#1}} }
\newcommand{\ta}[1]{\tau_{{}_{\!#1}} }
\newcommand{\n}[1]{\nu_{{}_{\!#1}} }
\newcommand{\til}[1]{\widetilde{#1}}

\newcommand{\link}{\mathop{\rm link}\nolimits}

\newcommand{\Hom}{\operatorname{Hom}}
\newcommand{\Tor}{\operatorname{Tor}}

\newcommand{\codim}{\operatorname{codim}}

\newcommand{\rank}{\operatorname{rank}}

\newcommand{\C}{\mathbb C}
\newcommand{\Q}{\mathbb Q}
\newcommand{\R}{\mathbb R}
\newcommand{\Z}{\mathbb Z}
\renewcommand{\k}{\mathbf k}
\renewcommand{\P}{\mathcal P}

\newcommand{\K}{\P}

\def\<{\langle}
\renewcommand{\>}{\rangle}
\renewcommand{\le}{\leqslant}
\renewcommand{\ge}{\geqslant}

\newcommand{\hatzero}{\hat 0}

\newcommand{\Me}{M_e}
\newcommand{\MG}{M_F}
\renewcommand{\H}{\widehat H}
\newcommand{\G}{F}

\begin{document}

\title{On the cohomology of torus manifolds}
\author{Mikiya Masuda}
\author{Taras Panov}\thanks{The second author was partially supported by
the Russian Foundation for Basic Research, grant no.~01-01-00546.}
\date{\today}
\address{Department of Mathematics, Osaka City University,
Sumiyoshi-ku, Osaka 558-8585, Japan}
\address{Department of Mathematics and Mechanics, Moscow State University,
Le\-nin\-s\-kiye Gory, Moscow 119992, Russia\newline
\emph{and }Institute for Theoretical and Experimental Physics, Moscow 117259,
Russia}
\email{masuda@sci.osaka-cu.ac.jp}
\email{tpanov@mech.math.msu.su}

\begin{abstract}
A \emph{torus manifold} is an even-dimensional manifold acted on by a
half-dimensional torus with non-empty fixed point set and some additional
orientation data. It may be considered as a far-reaching generalisation of
\emph{toric manifolds} from algebraic geometry. The orbit space of a torus
manifold has a rich combinatorial structure, e.g., it is a \emph{manifold
with corners} provided that the action is \emph{locally standard}. Here we
investigate relationships between the cohomological properties of torus
manifolds and the combinatorics of their orbit quotients. We show that the
cohomology ring of a torus manifold is generated by two-dimensional classes
if and only if the quotient is a \emph{homology polytope}. In this case we
retrieve the familiar picture from toric geometry: the equivariant
cohomology is the \emph{face ring} of the nerve simplicial complex and the
ordinary cohomology is obtained by factoring out certain linear forms. In a
more general situation, we show that the odd-degree cohomology of a torus
manifold vanishes if and only if the orbit space is \emph{face-acyclic}.
Although the cohomology is no longer generated in degree two under these
circumstances, the equivariant cohomology is still isomorphic to the face ring
of an appropriate \emph{simplicial poset}.
\end{abstract}

\maketitle

\section{Introduction}

Since the 1970s algebraic geometers have studied equivariant algebraic
compactifications of the \emph{algebraic torus} $(\C^*)^n$, nowadays known as
\emph{complete toric varieties}. The study quickly grew into a separate branch
of algebraic geometry, ``toric geometry", incorporating many topological and
convex-geometrical ideas and constructions, and producing a spectacular array
of applications.  A toric variety is a (normal) algebraic variety on which an
algebraic torus acts with a dense orbit. The variety and the action are fully
determined by a combinatorial object called a \emph{fan}~\cite{fult93}.

With the appearance of the pioneering work~\cite{da-ja91} of Davis and
Januszkiewicz in the beginning of the 1990s, the ideas of toric geometry have
started penetrating into topology. The orbit space of a non-singular
projective toric variety with respect to the
action of the compact torus $T^n\subset(\C^n)^*$ can be identified with the
simple polytope ``dual" to the corresponding fan. Moreover, the action of the
compact torus on a non-singular toric variety is ``locally standard", that
is, locally modelled by the standard action on $\C^n$. Davis and
Januszkiewicz took these two characteristic properties as a starting point
for their topological generalisation of toric varieties, namely \emph{quasitoric
manifolds}. A quasitoric manifold is a compact manifold $M^{2n}$ with a
locally standard action of $T^n$ whose orbit space is (combinatorially) a
simple polytope. (Davis and Januszkiewicz used the term ``toric manifold",
but by the time their work appeared the latter had already been used in
algebraic geometry as a synonym of ``non-singular toric variety".) According
to one of the main results of~\cite{da-ja91}, the cohomology ring of a
quasitoric manifold $M$ has the same structure as that of a non-singular
complete toric variety, and is isomorphic to the quotient of the
Stanley--Reisner face ring of the orbit space by certain linear forms. In
particular, the cohomology of $M$ is generated by degree-two elements.

In contrast, the convex-geometrical notion of polytope, while playing a
very important role in geometrical considerations related to toric geometry,
appears to be less relevant in the topological study of torus actions. The
orbit quotient $Q=M/T$ of a non-singular compact toric variety $M$
locally looks like the positive cone $\R^n_+$ and thereby acquires a
specific face decomposition. This combinatorial structure on $Q$ is known
to differential topologists as that of a \emph{manifold with corners}.
Moreover, all faces of $Q$, including $Q$ itself, and all their
intersections are acyclic. We call such a manifold with corners a
\emph{homology polytope}. It is a genuine polytope provided that the toric
variety is projective, but in general may fail to be so. This implies, in
particular, that the class of quasitoric manifolds does not include all
non-singular compact toric varieties (see~\cite[\S5.2]{bu-pa02} for more
discussion on the relationships between toric varieties and quasitoric
manifolds). On the other hand we might expect that all the
topological properties of quasitoric manifolds would still hold under a
weaker assumption that the orbit space of the torus action is a
homology polytope. This is justified by some results of the
present paper (see Theorem~\ref{hptpe}).

An alternative far-reaching topological generalisation of complete
non-\-sin\-gu\-lar toric varieties was introduced in~\cite{masu99}
and~\cite{ha-ma03} under the name of \emph{torus manifolds} (or \emph{unitary
toric manifolds} in the earlier terminology). A torus manifold is an
even-dimensional manifold $M$ acted on by a half-dimensional torus $T$ with
non-empty fixed point set; we also specify certain orientation data on $M$
from the beginning, in order to make certain isomorphisms canonical.
Particular examples of torus manifolds
include non-singular complete toric varieties (otherwise known as toric
manifolds) and the quasitoric manifolds of Davis and Januszkiewicz. On the other
hand, the conditions on the action are significantly weakened in comparison
to quasitoric manifolds. Surprisingly, torus manifolds admit a
combinatorial treatment similar to toric varieties. It relies on the notions
of \emph{multi-fans} and \emph{multi-polytopes}, developed in~\cite{ha-ma03}
as an alternative to fans associated with toric varieties.

The notion of torus manifold appears to be an appropriate concept for
investigating relationships between the topology of torus action and the
combinatorics of orbit quotient, which is the main theme of the current
paper. Our first main result (Theorem~\ref{hptpe}) measures the
extent of the analogy between the cohomological structure of non-singular
complete toric varieties and torus manifolds:
\begin{theorem*}
The cohomology of a torus manifold $M$ is generated by its degree-two part if
and only if $M$ is locally standard and the orbit space $Q$ is a homology
polytope.
\end{theorem*}
The cohomology ring itself may also be calculated and has a structure familiar
from toric geometry: it is isomorphic to the Stanley--Reisner face ring of
$Q$ modulo certain linear forms.

Next we study a more general class of torus manifolds: those with vanishing
odd-degree cohomology. Under these circumstances the equivariant cohomology
of $M$ is a free finitely generated module over the equivariant cohomology of
point, $H^*_T(pt)=\Z[t_1,\ldots,t_n]$. This condition is known to algebraists
as \emph{Cohen--\-Mac\-au\-lay\-ness} and is equivalent to $M$ being
\emph{equivariantly formal} in the terminology of~\cite{GS}. The orbit space
of a torus manifold with $H^{odd}(M)=0$ may fail to be a
homology polytope, as a simple example of torus acting on an even-dimensional
sphere shows (see Example~\ref{2nsphere}). We introduce a weaker notion
of \emph{face-acyclic} manifold with corners $Q$, in which all the faces are
still acyclic, but their intersections may fail to be connected, and prove
\begin{theorem*}
The odd-degree cohomology of $M$ vanishes if and only if $M$ is locally
standard and the orbit space $Q$ is face-acyclic.
\end{theorem*}
This result is stated as Theorem~\ref{cohfa} in our paper.  We also show that
the equivariant cohomology is isomorphic to the face ring of the simplicial
poset of faces of $Q$ and identify the ordinary cohomology accordingly
(Theorem~\ref{theo:face CM} and Corollary~\ref{coro:stcoh}). The face ring
of a simplicial poset is not generated by its degree-two elements in general.

At the end we prove Stanley's conjecture on the characterisation of
$h$-vectors of Gorenstein* simplicial posets in the particular case of face
posets of orbit quotients for torus manifolds (Theorem~\ref{theo:betti}).
Unlike the case of Gorenstein* simplicial complexes (which can be considered
as an ``algebraic approximation" to triangulations of spheres), the
conditions for an integer vector to be an $h$-vector of a Gorenstein*
simplicial poset are relatively weak. Such an $h$-vector must have
non-negative entries $h_i$ and satisfy the \emph{Dehn--Sommerville
equations} $h_i=h_{n-i}$, \ $i=0,\ldots,n$.
There are no other conditions for odd~$n$. In even dimensions there is one
other troublesome condition; the middle-dimensional entry of the $h$-vector
must be even if at least one other entry is zero. It is not hard to check
that these conditions are sufficient, by providing the corresponding examples
of simplicial posets. We show that these simplicial posets can be realised as
the face posets of orbit quotients for torus manifolds with $H^{odd}(M)=0$
(so that the $h$-vectors of posets are the even Betti vectors of torus
manifolds).  Stanley's conjecture~\cite{stan91} was that those three
conditions are also necessary. In this paper we establish the necessity for
$h$-vectors of posets associated to torus manifolds with $H^{odd}(M)=0$. This
is done through the calculation of the Stiefel--Whitney classes of torus
manifolds. Similar topological ideas were used by the first author to prove
the Stanley conjecture in full generality in~\cite{masu??}.

We note that the characterisation of $h$-vectors for Gorenstein* simplicial
complexes, as well as for sphere triangulations, remains wide open.

The paper is organised as follows. In Section~\ref{preli} we establish the
notation concerning torus actions on manifolds and prove three pivotal
statements (Lemmas \ref{lemm:freeness}--\ref{sconn}) describing different
properties of fixed point sets. In Section~\ref{sect:torus manifolds} we
introduce the concept of torus manifold, give a few examples, and establish
some basic facts about them. In Section~\ref{local} we discuss locally
standard torus actions. The main result here is Theorem~\ref{theo:local
standardness} showing that a torus manifold $M$ is locally standard provided
that $H^{odd}(M)=0$. We also introduce a canonical model for a torus manifold
with given orbit space $Q$ and the distribution of circle subgroups fixing
characteristic submanifolds. Then we show that a torus manifold is
equivariantly diffeomorphic to its canonical model provided that $H^2(Q)=0$.
This extends the corresponding result for quasitoric manifolds due to Davis
and Januszkiewicz. In Section~\ref{facer} we develop the necessary apparatus
of ``combinatorial commutative algebra". Here we introduce face rings of
manifolds with corners and simplicial posets, and list their main algebraic
properties. We try not to overload the notation with poset terminology, but a
reader familiar with posets will recognise the notions of (semi)lattice,
meet, join, etc. In Section~6 we turn to the equivariant cohomology of torus
manifolds. We introduce certain key concepts and construct a map from the
face ring of the orbit quotient to the equivariant cohomology of the torus
manifold, which is later shown to be an isomorphism under certain conditions.
Sections 7--9 contain the proofs of the main results quoted above.
In Section~\ref{gopos} we prove the above
mentioned particular case of Stanley's conjecture on Gorenstein* simplicial
posets.

\section{Preliminaries}\label{preli}

We start with recalling some basic theory of $G$-spaces, referring
to~\cite[Ch.~II]{bred72} for the proofs of the corresponding statements. Let
$X$ be a topological space with a left action of a compact topological group
$G$. The action is \emph{effective} if unit is the only element of $G$ that
acts trivially, and is \emph{free} if the \emph{isotropy subgroup}
$G_x=\{g\in G\colon gx=x\}$ is trivial for all $x\in X$. The fixed point set
is denoted $X^G$.  There exists a contractible free right $G$-space $EG$
called the \emph{universal $G$-space}; the quotient $BG:=EG/G$ is called the
\emph{classifying space} for free $G$-actions. The product $EG\times X$ is a
free left $G$-space by $g\cdot(e,x)=(eg^{-1},gx)$; the quotient $EG\times_G
X:=(EG\times X)/G$ is called the \emph{Borel construction} on $X$ or the
\emph{homotopy quotient} of~$X$. The \emph{equivariant cohomology} with
coefficients in a ring $\k$ is defined as
\[
  H^*_G(X;\k):=H^*(EG\times_G X;\k).
\]
The map $\rho$ collapsing $X$ to a point induces a homomorphism
\begin{equation}\label{rhost}
  \rho^*\colon H^*_G(pt;\k)=H^*(BG;\k)\to H^*_G(X;\k)
\end{equation}
thereby defining a canonical $H^*(BG;\k)$-module structure on $H^*_G(X;\k)$.
The
Borel construction can also be applied to a $G$-vector bundle. For instance,
if $E$ is an oriented $G$-vector bundle over a $G$-space $X$, then the Borel
construction on $E$ produces an oriented vector bundle over $EG\times_G X$
and its Euler class is called the \emph{equivariant Euler class} of $E$ and
denoted by $e^G(E)$.  Note that $e^G(E)$ lies in $H^{*}_G(X;\Z)$. Below we
use integer coefficients, unless another coefficient ring is specified.

If $G$ is a commutative group (e.g., a compact torus $T=T^k$), then the
notions of left and right $G$-spaces coincide.  As is well known, $H^*(BT)$
is a polynomial ring in $k$ variables of degree two, in particular
$H^{odd}(BT)=0$. All manifolds $M$ in this paper are closed connected smooth
and orientable.

\begin{lemma}\label{lemm:freeness}
Let $M$ be a manifold with a smooth action of $T$ such that the
fixed point set $M^T$ is finite and non-empty.  Then $H^*_T(M)$ is free as an
$H^*(BT)$-module if and only if $H^{odd}(M)=0$. In this case $H^*_T(M)\cong
H^*(BT)\otimes H^*(M)$ as $H^*(BT)$-modules.
\end{lemma}
\begin{proof}
Assume $H^{odd}(M)=0$. Then the Serre spectral sequence of the fibration
$ET\times_T M\to BT$ collapses and $H^*(M)$ has no torsion, so
$H^*_T(M)$ is isomorphic to $H^*(BT)\otimes H^*(M)$ and thus is a free
$H^*(BT)$-module.  This proves the ``if" part.

To prove the ``only if" part, we use the Eilenberg--Moore spectral sequence
of the bundle $ET\times_T M\to BT$ with fibre $M$.  It converges to $H^*(M)$
and has
$$
  E_2^{*,*}=\Tor^{*,*}_{H^*(BT)}\bigl(H^*_T(M),\Z\bigr).
$$
Since $H_T^*(M)$ is free as an $H^*(BT)$-module,
we have
\begin{align*}
  \Tor^{*,*}_{H^*(BT)}\bigl(H^*_T(M),\Z\bigr)&=
  \Tor^{0,*}_{H^*(BT)}\bigl(H^*_T(M),\Z\bigr)\\
  &=H^*_T(M)\otimes_{H^*(BT)}\Z\\
  &=  H^*_T(M)\bigr/(\l^*(H^{>0}(BT))).
\end{align*}
Therefore, $E_2^{0,*}=H^*_T(M)\bigr/(\l^*(H^{>0}(BT)))$ and $E_2^{-p,*}=0$
for $p>0$. It follows that the Eilenberg--Moore spectral sequence collapses
at the $E_2$ term and
\begin{equation} \label{eqn:H(M)}
  H^*(M)=H^*_T(M)/(\l^*(H^{>0}(BT))).
\end{equation}
On the other hand, it follows from the localisation theorem (see
\cite{Hsiang}) that the kernel of the restriction map
\[
  H^*_T(M)\to H^*_T(M^T)=H^*(BT)\otimes H^*(M^T)
\]
is the $H^*(BT)$-torsion subgroup and hence the restriction map is injective
in our case. Therefore $H^{odd}_T(M)=0$ because $M^T$ is a finite set of
isolated points. This fact together with (\ref{eqn:H(M)}) proves that
$H^{odd}(M)=0$.
\end{proof}

Two classes of $T$-manifolds, namely those having zero odd degree cohomology
or even cohomology generated in degree two, are of particular importance in
this paper. Next we prove two technical lemmas showing that these
cohomological properties are inherited by the fixed point set $M^H$ for any
subtorus $H\subseteq T$. These lemmas will be used in inductive arguments
later in the paper.

\begin{lemma}\label{lemm:odd=0}
Let $M$ be a $T$-manifold, $H$ a subtorus of $T$ and $N$ a connected
component of $M^H$. If $H^{odd}(M)=0$, then $H^{odd}(N)=0$ and
$N^T\not=\varnothing$.
\end{lemma}
\begin{proof}
We first prove that $H^{odd}(M^H)=0$. Note that for a generic circle subgroup
$S\subseteq H$ we have $M^S=M^H$. Let $p$ be a prime and $G$ be an order $p$
subgroup in $S$. The induced action of $G$ on $H^*(M)$ is trivial because $G$
is contained in the connected group $S$. Then $\dim H^{odd}(M^G;\Z/p)\le\dim
H^{odd}(M;\Z/p)$ by~\cite[Theorem~VII.2.2]{bred72}. Therefore,
$H^{odd}(M^G;\Z/p)=0$ by the assumption.  Repeating the same argument for
$M^G$ with the induced action of $S/G$, which is again a circle group, we
conclude that $H^{odd}(M^G;\Z/p)=0$ for any $p$-subgroup $G$ of $S$. However,
$M^G=M^{S}=M^H$ if the order of $G$ is sufficiently large, so we have
$H^{odd}(M^H;\Z/p)=0$. Since $p$ is an arbitrary prime, this implies that
$H^{odd}(M^H)=0$.

Now since $H^{odd}(N)=0$, the Euler characteristic $\chi(N)$ of $N$ is
non-zero. As is well-known $\chi(N)=\chi(N^T)$, which implies that $N^T$ is
non-empty.
\end{proof}

\begin{lemma}\label{sconn}
Let $M,H,N$ be as in Lemma~\ref{lemm:odd=0}.  If $H^*(M)$ is generated by its
degree-two part (as a ring), then the restriction map $H^*(M)\to H^*(N)$ is
surjective; in particular, $H^*(N)$ is also generated by its degree-two part.
\end{lemma}
\begin{proof}
Since $H^{odd}(M)=0$, we have $H^{odd}(N)=0$ by Lemma~\ref{lemm:odd=0}; so it
suffices to prove that the restriction map $H^*(M;\Z/p)\to H^*(N;\Z/p)$ is
surjective for any prime $p$.

The argument below is similar to that used in the proof of Theorem~VII.3.1
in~\cite{bred72}.
As in the proof of Lemma~\ref{lemm:odd=0}, let $S$ be a generic circle
subgroup of $H$ (so that $M^S=M^H$) and let $G$ be the subgroup of $S$
of prime order $p$. Then the restriction map $H^k_G(M;\Z/p) \to
H^k_G(M^G;\Z/p)$ is an isomorphism for sufficiently large~$k$
by~\cite[Theorem~VII.1.5]{bred72}.
Hence, for any connected component $N'$ of $M^G$ the restriction $r\colon
H^k_G(M;\Z/p) \to H^k_G(N';\Z/p)$ is surjective if $k$ is sufficiently large.
Now consider the commutative diagram
\[
\begin{CD}
  H^*_{G}(M;\Z/p) @>r>> H^*_{G}(N';\Z/p)@.\;\cong H^*(BG;\Z/p)
  \otimes H^*(N';\Z/p)\\
  @VVV @VVV @.\\
  H^*(M;\Z/p) @>s>> H^*(N';\Z/p)@.
\end{CD}.
\]
Choose a basis $v_1,\ldots,v_d\in H^2(M;\Z/p)$; then these elements are
multiplicative generators for $H^*(M;\Z/p)$. Since
$H^{odd}(M;\Z/p)=H^{odd}(M^G;\Z/p)=0$ and $\chi(M)=\chi(M^T)=\chi(M^G)$, we
have $\sum\dim H^i(M;\Z/p)=\sum\dim H^i(M^G;\Z/p)$.
By~\cite[Theorem~VII.1.6]{bred72} the Serre spectral sequence of the
fibration $EG\times_G M\to BG$ collapses. Therefore, the vertical map
$H^*_G(M;\Z/p)\to H^*(M;\Z/p)$ in the above diagram is surjective. Let
$\xi_j\in H^*_{G}(M;\Z/p)$ be a lift of $v_j$, and $w_j:=s(v_j)$. Let $t$ be
a generator of $H^2(BG;\Z/p)\cong\Z/p$.  Since the above diagram is
commutative and $H^1(N';\Z/p)=0$ by Lemma~\ref{lemm:odd=0}, we have
$r(\xi_j)=\alpha_jt+w_j$ for some $\alpha_j\in \Z/p$. Now let $a\in
H^*(N';\Z/p)$ be an arbitrary element. Then there exist $\ell$ and a
polynomial $P(\xi_1,\dots,\xi_d)$ such that
\[
  r\bigl(P(\xi_1,\dots,\xi_d)\bigr)=t^{\ell}a.
\]
On the other hand,
\[
  r\bigl(P(\xi_1,\dots,\xi_d)\bigr)=P(\alpha_1t+w_1,\dots,\alpha_dt+w_d)=
  \sum_{k\ge0}t^kQ_k(w_1,\dots,w_d)
\]
for some polynomials $Q_k$. Therefore, $a=Q_\ell(w_1,\dots,w_d)$, the
restriction map $H^*(M;\Z/p)\to H^*(N';\Z/p)$ is surjective, and
$H^*(N';\Z/p)$ is generated by the degree-two elements $w_1,\dots,w_d$.

Now we can repeat the same argument for $N'$ with the induced action of
$S/G$, which is again a circle group. It follows that the restriction map
$H^*(M;\Z/p)\to H^*(N';\Z/p)$ is surjective for any connected component $N'$
of $M^G$ with $G$ any $p$-subgroup of $S$. However, if the order of $G$ is
sufficiently large, then $M^G=M^S=M^H$ and hence $N'=N$, so it follows
that the restriction map $H^*(M;\Z/p)\to H^*(N;\Z/p)$ is surjective for any
connected component $N$ of $M^{H}$. Since the prime $p$ is arbitrary, the
proof is finished.
\end{proof}

\section{Torus manifolds}\label{sect:torus manifolds}

The notion of torus manifold was introduced in
\cite{ha-ma03} and \cite{masu99}, and here we follow the notation of these
papers with some additional specifications.

A \emph{torus manifold} is a $2n$-dimensional closed connected orientable
smooth manifold $M$ with an effective smooth action of an $n$-dimensional
torus $T=(S^1)^n$ such that $M^T\not=\varnothing$.  Since $\dim M=2\dim T$
and $M$ is compact, the fixed point set $M^T$ is a finite set of isolated
points.

A codimension-two connected component of the set fixed pointwise by a
circle subgroup of $T$ is called a \emph{characteristic submanifold}
of~$M$.  The existence of a $T$-fixed
point is required for the definition of characteristic submanifold
in \cite{ha-ma03} and \cite{masu99} but not in this paper.
However, when $H^{odd}(M)=0$, these two definitions agree by
Lemma~\ref{lemm:odd=0}.

Since $M$ is compact, there are only finitely many characteristic
submanifolds, and we denote them by $M_i$, \ $i=1,\dots,m$.
Each characteristic submanifold $M_i$ is orientable as a connected
component of the fixed
point set for a circle action on an orientable manifold.
Following~\cite{bu-ra01}, we say that $M$ is \emph{omnioriented} if an
orientation is specified for $M$ and for every characteristic submanifold
$M_i$. There are $2^{m+1}$ choices of omniorientations. It is extremely
convenient, although not absolutely necessary to assume that all torus
manifolds are omnioriented (in~\cite{ha-ma03} a choice of omniorientation
for characteristic submanifolds was
a part of the definition of torus manifold).

Here are two typical examples of torus manifolds.

\begin{example} \label{CPn}
A complex projective space $\C P^n$ has a natural $T$-action defined in the
homogeneous coordinates by
\[
  (t_1,\dots,t_n)\cdot(z_0:z_1:\dots:z_n)=(z_0:t_1z_1:\dots:t_nz_n).
\]
It has $(n+1)$ characteristic submanifolds $\{z_0=0\},\ldots,\{z_n=0\}$ and
$(n+1)$ fixed points $(1:0:\dots:0),\ldots,(0:\dots:0:1)$. In this example
the intersection of any set of characteristic submanifolds is connected.
\end{example}

\begin{example}\label{2nsphere}
Let $S^{2n}$ be the $2n$-sphere identified with the following subset in
$\C^{n}\times\R$:
$$
  \bigl\{ (z_1,\dots,z_n,y)\in\C^n\times\R\colon|z_1|^2+\dots+|z_n|^2+y^2=1
  \bigr\}.
$$
Define a $T$-action by
$$
  (t_1,\dots,t_n)\cdot(z_1,\dots,z_n,y)=(t_1z_1,\dots,t_nz_n,y).
$$
It has $n$ characteristic submanifolds $\{z_1=0\},\ldots,\{z_n=0\}$, and two
fixed points $(0,\dots,0,\pm 1)$. The intersection of any $k$ characteristic
submanifolds is connected if $k\le n-1$, but consists of two disjoint
fixed points if $k=n$.
\end{example}

If $M$ is an (omnioriented) torus manifold, then both $M$ and $M_i$ are
oriented, and the Gysin homomorphism $H_T^*(M_i)\to H_T^{*+2}(M)$ in
equivariant cohomology is defined. Denote by $\tau_i\in H^2_T(M)$ the image
of the identity element in $H^0_T(M_i)$. We may think of $\tau_i$ as the
Poincar\'e dual of $M_i$ in equivariant cohomology.

\begin{proposition}[See section 1 of \cite{masu99}]\label{prop:a_i}
Let $M$ be a torus manifold.
\begin{itemize}

\item[1.] For each characteristic submanifold $M_i$ with
$(M_i)^T\not=\varnothing$,
there is a unique element $a_i\in H_2(BT)$ such that
$$
  \l^*(t)=\sum_i\<t,a_i\>\tau_i\quad
  \text{modulo $H^*(BT)$-torsions}
$$
for any element $t\in H^2(BT)$.  Here the sum is taken over all
characteristic submanifolds $M_i$ with $(M_i)^T\not=\varnothing$ and
$\rho^*$ denotes the homomorphism~\eqref{rhost}.

\item[2.] The circle subgroup fixing $M_i$ with $(M_i)^T\not=\varnothing$
coincides with the one determined by $a_i\in H_2(BT)$ through the
identification $H_2(BT)=\Hom(S^1,T)$.

\item[3.] If $n$ different characteristic submanifolds
$M_{i_1},\dots,M_{i_n}$ have a $T$-fixed point in their intersection, then
the elements $a_{i_1},\ldots,a_{i_n}$ form a basis of $H_2(BT)$ over $\Z$.
\end{itemize}
\end{proposition}

The next lemma provides a sufficient cohomological condition for the
intersections of characteristic submanifolds to be connected (compare
Examples~\ref{CPn} and~\ref{2nsphere}).

\begin{lemma}\label{cmcon}
Suppose that $H^*(M)$ is generated in degree two. Then all non-empty multiple
intersections of the characteristic submanifolds are connected and have
cohomology generated in degree two.
\end{lemma}
\begin{proof}
Since every characteristic submanifold $M_i$ is a connected component of the
fixed point set of a circle subgroup of $T$, the cohomology $H^*(M_i)$ is
generated by the degree-two part and the restriction map $H^*(M)\to H^*(M_i)$
is onto by Lemma~\ref{sconn}. It follows that the restriction map $H^*_T(M)\to
H^*_T(M_i)$ in equivariant cohomology is also onto.

Now we prove the connectedness of multiple intersections.  Suppose that
$M_{i_1}\cap\dots\cap M_{i_k}\ne\varnothing$, \ $(1<k\le n)$, and pick a
connected component $N$ of the intersection.  Since $N$ is fixed by a
subtorus, it contains a $T$-fixed point by Lemma~\ref{lemm:odd=0}.
For each $i\in\{i_1,\dots,i_k\}$ there are
embeddings $\f_i\colon N\to M_i$, $\psi_i\colon M_i\to M$, and the
corresponding Gysin homomorphisms in equivariant cohomology:
\[
\begin{CD}
  H^0_T(N) @>\f_{i_!}>> H^{2k-2}_T(M_i) @>\psi_{i_!}>> H^{2k}_T(M).
\end{CD}
\]
Since the restriction
$\psi_i^*\colon H^*_T(M)\to H^*_T(M_i)$ is surjective, we have
$\f_{i_!}(1)=\psi^*_i(u)$ for some $u\in H^{2k-2}_T(M)$.  Now we calculate
\[
  (\psi_i\circ\f_i)_!(1)=\psi_{i_!}(\f_{i_!}(1))=
  \psi_{i_!}\bigl(\psi_i^*(u)\bigr)=\psi_{i_!}(1)u=\tau_i u.
\]
Hence, $(\psi_i\circ\f_i)_!(1)$ is divisible by $\tau_i$ for every
$i\in\{i_1,\dots,i_k\}$. By Proposition 3.4 of \cite{masu99}, the degree-$2k$
part of $H^*_T(M)$ is additively generated by the monomials
$\tau^{k_1}_{j_1}\dots \tau^{k_p}_{j_p}$ such that $M_{j_1}\cap\dots\cap
M_{j_p}\not=\varnothing$ and $k_1+\dots+k_p=k$. It follows that
$(\psi_i\circ\f_i)_!(1)$ is a non-zero integral multiple of $\tau_{i_1}\dots
\tau_{i_k}\in H^{2k}_T(M)$. By the definition of Gysin map,
$(\psi_i\circ\f_i)_!(1)$ goes to zero under the restriction map $H_T^*(M)\to
H_T^*(x)$ for every point $x\in (M\backslash N)^T$. On the other hand, the
image of $\tau_{i_1}\dots \tau_{i_k}$ under the restriction map $H_T^*(M)\to
H_T^*(x)$ is non-zero for every $T$-fixed point $x\in M_{i_1}\cap\dots\cap
M_{i_k}$. Thus, $N$ is the only connected component of the latter
intersection. The fact that $H^*(N)$ is generated by its degree-two part
follows from Lemma~\ref{sconn}.
\end{proof}

\section{Locally standard torus manifolds and orbit spaces}\label{local}
\subsection{Locally standardness}

We say that a torus manifold $M$ is \emph{locally standard} if every point in
$M$ has an invariant neighbourhood $U$ weakly equivariantly diffeomorphic to
an open subset $W\subset\C^n$ invariant under the standard $T^n$-action on
$\C^n$. The latter means that there is an automorphism $\psi\colon T\to T$
and a diffeomorphism $f\colon U\to W$ such that $f(ty)=\psi(t)f(y)$ for all
$t\in T$, \ $y\in U$.

The following statement gives a sufficient cohomological condition for
local standardness.
\begin{theorem}\label{theo:local standardness}
A torus manifold $M$ with $H^{odd}(M)=0$ is locally standard.
\end{theorem}
\begin{proof}
We first show that there are no non-trivial finite isotropy subgroups for the
$T$-action on~$M$. Assume the opposite, i.e., the isotropy group $T_x$ is
finite and non-trivial for some $x\in M$. Then $T_x$ contains a non-trivial
cyclic subgroup $G$ of some prime order~$p$. Let $N$ be the connected component
of $M^G$ containing $x$.  Since $N$ contains $x$ and $T_x$ is finite, the
principal isotropy group of $N$ is finite. Like in the proof of
Lemma~\ref{lemm:odd=0}, it follows from~\cite[Theorem~VII.2.2]{bred72} that
$H^{odd}(N;\Z/p)=0$. In particular, the Euler characteristic of $N$ is
non-zero, and therefore, $N$ has a $T$-fixed point, say $y$. The tangential
$T$-representation $\mathcal T_y M$ at $y$ is faithful, $\dim M=2\dim T$ and
$\mathcal T_y N$ is a proper $T$-subrepresentation of $\mathcal T_yM$. It
follows that there is a subtorus $T'$ (of positive dimension) which fixes
$\mathcal T_yN$ and does not fix the complement of $\mathcal T_yN$ in
$\mathcal T_yM$. Clearly, $T'$ is the principal isotropy group of $N$, which
contradicts the above observation that the principal isotropy group of $N$ is
finite.

If the isotropy group $T_x$ is trivial, $M$ is obviously locally standard
near~$x$. Suppose that $T_x$ is non-trivial. Then it cannot be finite and
therefore, $\dim T_x>0$.  Let $H$ be the identity component of $T_x$, and $N$
the connected component of $M^{H}$ containing $x$. By Lemma~\ref{lemm:odd=0},
$N$ has a $T$-fixed point, say $y$. Looking at the tangential representation
at $y$, we observe that the induced action of $T/H$ on $N$ is effective. By
the previous argument, no point of $N$ has a non-trivial finite isotropy
group for the induced action of $T/H$, which implies that $T_x=H$. Since $x$
and $y$ are both in the same connected component $N$ fixed pointwise by
$T_x$, the $T_x$-representation in $\mathcal T_x M$ agrees with the
restriction of the tangential $T$-representation $\mathcal T_yM$ to $T_x$.
This implies that $M$ is locally standard near~$x$.
\end{proof}

In the rest of this section we assume that $M$ is locally standard.

Let $Q:=M/T$ denote the orbit space of $M$ and $\pi\colon M\to Q$ the
quotient projection.
Since $M$ is locally standard, any point in the orbit space $Q$
has a neighbourhood diffeomorphic to an open subset
in the positive cone
$$
  \R^n_\ge=\{(y_1,\dots,y_n)\in\R^n\colon y_i\ge0,\ i=1,\dots,n\}.
$$
This identifies $Q$ as a \emph{manifold with corners}, see
e.g.~\cite[\S6]{davi83}, and faces of $Q$ can be defined in a natural way.
The vertices of $Q$, that is, the $0$-dimensional faces, correspond to
the $T$-fixed points of $M$ through the quotient projection $\pi$.
Codimension one faces of $Q$ are called the \emph{facets} of $Q$.
They are the $\pi$ images of characteristic submanifolds $M_i$,
\ $i=1,\dots,m$.  We set $\F_i:=\pi(M_i)$.
We refer to a non-empty intersection of $k$ facets as a codimension-$k$
\emph{preface}, \ $k=1,\dots,n$.
In general, prefaces of codimension $>1$
may fail to be connected (see Example~\ref{2nsphere}).
Faces are connected components of prefaces.
We also regard $Q$ itself as a
codimension-zero face; other faces are called \emph{proper faces}.
If $H^{odd}(M)=0$, then every face has a vertex by Lemma~\ref{lemm:odd=0}.
Moreover, if $H^*(M)$ is generated in degree two, then all prefaces are
connected by Lemma~\ref{cmcon}; so prefaces are faces in this case.

A space
$X$ is \emph{acyclic} if $\widetilde H_i(X)=0$ for all $i$. We say that $Q$
is \emph{face-acyclic} if all of its faces (including $Q$ itself) are acyclic.
It is not difficult to see that if $Q$ is face-acyclic,
then every face of $Q$ has a vertex.
We call $Q$ a \emph{homology polytope} if all its prefaces are
acyclic (in particular, connected), in other words, $Q$ is a homology polytope
if and only if it is face-acyclic and all non-empty multiple intersections of
characteristic submanifolds are connected.

\begin{remark}
A \emph{simple convex polytope} is an example of a manifold with corners and
is a homology polytope. A \emph{quasitoric manifold}~\cite{da-ja91},
\cite{bu-pa02} can be defined as a locally standard torus manifold whose
orbit space is a simple convex polytope with the standard face structure.
\end{remark}

\begin{example}\label{CPn orbit}
Torus manifold $\C P^n$ with the $T$-action from Example~\ref{CPn} is
locally standard and the map
\[
  (z_0:z_1:\dots:z_n) \to
  \frac{1}{\sum_{i=0}^n |z_i|^2}(|z_1|^2,\dots,|z_n|^2)
\]
induces a face preserving homeomorphism from the orbit space $\C P^n/T$ to a
standard $n$-simplex. The latter is a simple polytope, in particular, a
homology polytope.
\end{example}

\begin{example}\label{S2n orbit}
Torus manifold $S^{2n}$ with the $T$-action from Example~\ref{2nsphere}
is also locally standard and the map
\[
  (z_1,\dots,z_n,y)\to (|z_1|,\dots,|z_n|,y)
\]
induces a face preserving homeomorphism from the orbit space
$S^{2n}/T$ to the space
\[
  \{ (x_1,\dots,x_n,y)\in \R^{n+1}\colon x_1^2+\dots+x_n^2+y^2=1,\ x_1\ge 0,
  \dots,x_n\ge 0\}.
\]
This space is not a homology polytope, but is a face-acyclic manifold with
corners.
\end{example}

\subsection{Canonical model}

In this paragraph we reconstruct the torus manifold $M$ from the orbit space
$Q$ and a map $\Lambda$ defined below using a ``canonical model"
$M_Q(\Lambda)$, which generalises a result of
Davis--Januszkiewicz~\cite[Prop.~1.8]{da-ja91}.

Remember that $M_i=\pi^{-1}(\F_i)$ is fixed by a circle
subgroup of $T$.  We choose a map
\begin{equation}
\label{eqn:characteristic map}
  \L\colon \{\F_1,\dots,\F_m\} \to H_2(BT)=\Hom(S^1,T)\cong \Z^n
\end{equation}
such that $\Lambda(\F_i)$ is primitive and determines the circle subgroup
of $T$ fixing $M_i$.  When $M_i$ has a $T$-fixed point, $\Lambda(\F_i)$
coincides with the element $a_i$ introduced in Proposition~\ref{prop:a_i} up
to sign. The following lemma follows immediately from the local standardness
of $M$.

\begin{lemma}\label{nonsi}
If $\F_{i_1}\cap \dots\cap \F_{i_k}$ is non-empty, then
$\L(\F_{i_1}),\dots,\L(\F_{i_k})$ is a part of basis for the integral lattice
$\Hom(S^1,T)\cong\Z^n$.
\end{lemma}

Given a point $x\in Q$, the smallest face which contains $x$ is an
intersection $\F_{i_1}\cap\dots\cap \F_{i_k}$ of some facets, and we define
$T(x)$ to be the subtorus of $T$ generated by the circle subgroups
corresponding to $\L(\F_{i_1}),\dots,\L(\F_{i_k})$.
Now introduce the identification space
\begin{equation}\label{idspa}
  M_Q(\L):=T\times Q/\!\sim,
\end{equation}
where $(t,x)\sim(t',x')$ if and only if $x=x'$ and $t^{-1}t'\in T(x)$. The
space $M_Q(\L)$ admits a natural action of $T$ and is a closed manifold
(this follows from Lemma~\ref{nonsi} and the fact that $Q$ is a manifold
with corners).
The following is a straightforward generalisation of
a~\cite[Prop.~1.8]{da-ja91}.

\begin{lemma}\label{acydj}
Let $M$ be a locally standard torus manifold with orbit space $Q$, and
the map $\L$ defined by~\eqref{eqn:characteristic map}. If $H^2(Q)=0$,
then there is an equivariant homeomorphism
\[
  M_Q(\L)\to M
\]
covering the identity on~$Q$.
\end{lemma}
\begin{proof}
The idea is to construct a continuous map $f\colon T\times Q\to M$ taking
$T\times q$ onto $\pi^{-1}(q)$ for each point $q\in Q$. This is done by
subsequent ``blowing up the singular strata". The condition on the second
cohomology group guarantees that the resulting principal $T$-bundle over $Q$
is trivial. Then the map $f$ descends to the required equivariant
homeomorphism. See~\cite{da-ja91} for details.
\end{proof}

\begin{remark}
Like in the case of quasitoric manifolds, it follows that a torus manifold
whose orbit quotient $Q$ satisfies $H^2(Q)=0$ is determined by $Q$ and~$\L$.
\end{remark}

\section{Face rings of manifolds with corners and simplicial posets}
\label{facer}

Before we proceed with describing the ordinary and equivariant cohomology
rings of torus manifolds we need an algebraic digression. Here we review a
notion of face ring generalising the classical Stanley--Reisner face
ring~\cite{stan96} to combinatorial structures more general than simplicial
complexes.  We consider two cases, which are in a sense dual to each other:
\lq\lq nice" manifolds with corners and simplicial posets. The latter one is
more general, however the former one is more convenient for applications to
torus manifolds. The face ring of a manifold with corners is also a little
easier to visualise, so we start with considering this case.

The relationship between nice manifolds with corners and simplicial posets is
similar to that between simple polytopes and simplicial complexes.  Face
rings of simplicial posets were introduced and studied in~\cite{stan91}.
Most of the statements in this section follow from the general theory of
ASL's (\emph{algebras with straightening law}) and \emph{Hodge algebras} as
explained in~\cite{stan91} and~\cite[Ch.~7]{br-he98}, however our treatment
is independent and geometrical.

\subsection{Nice manifolds with corners}

To begin, we assume that $Q$ is a homology polytope (or even a simple convex
polytope) with $m$ facets $\F_1,\dots,\F_m$.  Let $\k$ be a ground
commutative ring with unit, and assign a degree-two polynomial
generator~$\v{Q_i}$ to each facet~$Q_i$. We refer to the quotient ring
\[
  \k[Q]=
  \k\bigl[\v{\F_1},\dots,\v{\F_m}\bigr]\bigl/\bigl(\v{\F_{i_1}}\!\!\cdots
  \v{\F_{i_k}}=0
  \quad\text{if}\quad \F_{i_1}\cap\dots\cap \F_{i_k}=\varnothing\bigr).
\]
as the \emph{face ring} of $Q$. In coincides with the Stanley--Reisner face
ring~\cite{stan96} of the nerve simplicial complex~$K$.

For arbitrary pair of faces $G,H$ of $Q$ the intersection $G\cap H$ is a
unique maximal face contained in both $G$ and $H$. There is also a unique
minimal face that contains both $G$ and $H$, which we denote $G\vee H$. Let
$\k[\v{\G}\colon \G\text{ a face}]$ be the graded polynomial ring with one
$2k$-dimensional generator $\v{\G}$ for every proper codimension-$k$ face
$\G$. We also identify $\v{Q}$ with the unit and $v_\varnothing$ with zero.
The following proposition gives another presentation of $\k[Q]$, by extending
both the set of generators and relations. It will be used for a subsequent
generalisation of $\k[Q]$ to arbitrary manifolds with corners.

\begin{proposition}\label{frpol}
There is a canonical isomorphism of rings
$$
  \k[\v{\G}\colon \G\text{ a face}]/\mathcal I_Q\cong\k[Q],
$$
where $\mathcal I_Q$ is the ideal generated by all elements
$$
  \v{G}\v{H}-\v{G\vee H}\v{G\cap H}.
$$
\end{proposition}
\begin{proof}
The identification is established by the map sending $\v{\G}$ to
$\prod_{\F_i\supseteq \G}\v{\F_i}$.
\end{proof}

Now let $Q$ be an arbitrary connected manifold with corners. We also assume
that $Q$ is \emph{nice}, that is, every codimension-$k$ face is contained in
exactly $k$ facets. Note that the orbit space of a locally standard torus
manifold is always nice. In a nice manifold with corners, all faces
containing a given face form a Boolean lattice (like in the case of
$\R^n_\ge$).

\begin{remark}
By the definition of manifold with corners, every codimension-$k$ face is
contained in at most $k$ facets. A 2-disc with one 0-face and one 1-face on
the boundary gives an example of manifold with corners which is not nice.
\end{remark}

The intersection of two faces $G$ and $H$ in a manifold with corners may be
disconnected, but every its connected component is a face of codimension
$\codim G+\codim H$. We regard $G\cap H$ as the set of its connected
components; so the notation $E\in G\cap H$ is used below for connected
components $E$ of the intersection.

\begin{proposition}
For every two faces $G$ and $H$ with non-empty intersection, there is a
unique minimal face $G\vee H$ that contains both $G$ and $H$.
\end{proposition}
\begin{proof}
Take any $E\in G\cap H$. The statement follows from the fact that
the poset of faces containing $E$ is a Boolean lattice.
\end{proof}

Now we use the interpretation from Proposition~\ref{frpol} to
introduce a more general version of $\k[Q]$.

\begin{definition}\label{frgen}
The \emph{face ring} $\k[Q]$
of a nice manifold with corners $Q$ is a graded ring defined by
$$
  \k[Q]:=\k[\v{\G}\colon \G\text{ a face}]/\mathcal I_Q,
$$
where $\deg \v{\G}=2\codim\G$ and
$\mathcal I_Q$ is the ideal generated by all elements
$$
  \v{G}\v{H}-
    \v{G\vee H}\cdot\!\!\!\sum_{E\in{G\cap H}}\!\!\!\v{E}.
$$

\end{definition}

If $G$ and $H$ are transversal, that is, $\codim
G\cap H=\codim G+\codim H$, then $G\vee H=Q$, so in $\k[Q]$ we get the
identity
$$
  \v{G}\v{H}=\sum_{E\in G\cap H}\v{E}.
$$

Below we give a sequence of statements describing algebraic properties of
$\k[Q]$ and emphasising its analogy with the classical Stanley--Reisner face
ring.

\begin{lemma}\label{order}
Every element $a\in\k[Q]$ can be written as a linear combination
\[
  a=\mathop{\sum_{G_1\supset\dots\supset G_q}}
  \limits_{\alpha_1,\ldots,\alpha_q}\!
  A(G_1\supset\dots\supset G_q;\alpha_1,\ldots,\alpha_q)\,
  \v{G_1}^{\alpha_1}\cdots\v{G_q}^{\alpha_q}
\]
with coefficients $A(G_1\supset\dots\supset
G_q;\alpha_1,\ldots,\alpha_q)\in\k$.  Here $\codim G_i=i$ and
$G_q$ is an inclusion minimal face, and the sum is taken over all
chains of faces $G_1\supset\cdots\supset G_q$ with all
non-negative integers $\alpha_i$.
\end{lemma}
\begin{proof}
We may assume that $a=\v{H_1}\v{H_2}\!\!\cdots\v{H_k}$ (some $H_i$ may
coincide), and it is enough to show that it can be written as
$\sum\v{G_1}\!\!\cdots\v{G_l}$ with $G_1\supseteq\cdots\supseteq G_l$ for
every summand (without making any assumptions on codimensions, but allowing
some $G_i$ to coincide). By induction we may assume that
$H_2\supseteq\cdots\supseteq H_k$. Now we apply the relation from
Definition~\ref{frgen} and replace $a$ by
\[
  \v{H_1\vee
  H_2}\Bigl(\sum_{E\in{H_1\cap H_2}}\!\!\!\v{E}\Bigr)
  \v{H_3}\!\!\cdots\v{H_k}.
\]
The first two faces in every summand above are ordered. Then we replace each
$\v{E}\v{H_3}$ by $\v{E\vee H_3}(\sum_{G\in{E\cap H_3}}\v{G})$. Since
$H_1\vee H_2\supseteq E\vee H_3$, we get the first three faces in a linear
order. Proceeding in this fashion we finally end up in a sum of monomials
corresponding to ordered sets of faces.
\end{proof}

We refer to the presentation from Lemma~\ref{order} as the \emph{chain
decomposition} of an element $a\in\k[Q]$.

For any vertex (0-face) $p\in Q$ we define the \emph{restriction map $s_p$}
by
$$
  s_p\colon\k[Q]\to\k[Q]/(\v{\G}\colon \G\not\ni p).
$$
The next observation is straightforward.

\begin{proposition}\label{restr}
The image $s_p(\k[Q])$ of the restriction map can be identified with the
polynomial ring $\k[\v{\F_{i_1}},\dots,\v{\F_{i_n}}]$ on $n$ degree-two
generators, where $\F_{i_1},\dots,\F_{i_n}$ are the $n$ different facets
containing~$p$.
\end{proposition}

\begin{lemma}\label{alres}
If every face of $Q$ has a vertex, then
the sum $s=\oplus_ps_p$ of restriction maps over all vertices
$p\in Q$ is a monomorphism from $\k[Q]$ to the sum of
polynomial rings.
\end{lemma}
\begin{proof}
Take a non-zero $a\in\k[Q]$ and write it as in Lemma~\ref{order}.  Fix a
monomial $\v{G_1}^{\alpha_1}\cdots\v{G_n}^{\alpha_n}$ entering the chain
decomposition with a non-zero coefficient, and consider the restriction $s_p$
to the vertex $p=G_n$. We claim that $s_p(a)\ne0$. Identify $s_p(\k[Q])$ with
the polynomial ring $\k[t_1,\dots,t_n]$ (so that $t_j:=\v{\F_{i_j}}$ in the
notation of Proposition~\ref{restr}). Then $s_p(\v{G_n})=t_1\cdots t_n$ and
we may also assume that $s_p(\v{G_j})=t_1\cdots t_j$, \ $j=1,\dots,n$. Hence,
$$
  s_p\bigl( \v{G_1}^{\alpha_1}\cdots\v{G_n}^{\alpha_n} \bigr)=
  t_1^{\alpha_1}(t_1t_2)^{\alpha_2}\cdots(t_1\cdots t_n)^{\alpha_n}.
$$
It follows that $s_p(a)\ne0$ unless some other monomial
$\v{H_1}^{\beta_1}\cdots\v{H_n}^{\beta_n}$ hits the same monomial in
$\k[t_1,\dots,t_n]$. Note that
$$
  s_p(\v{H_1}^{\beta_1}\cdots\v{H_n}^{\beta_n})=0\qquad\text{unless
  } H_k\supseteq G_n \text{ for }\beta_k\ne0 .
$$
Suppose
\begin{equation}\label{rpvgh}
  s_p\bigl( \v{G_1}^{\alpha_1}\cdots\v{G_n}^{\alpha_n} \bigr)=
  s_p\bigl( \v{H_1}^{\beta_1}\cdots\v{H_n}^{\beta_n} \bigr).
\end{equation}
We want to prove that $\v{G_1}^{\alpha_1}\cdots\v{G_n}^{\alpha_n}=
\v{H_1}^{\beta_1}\cdots\v{H_n}^{\beta_n}$, that is, $\alpha_i=\beta_i$ and
$G_i=H_i$ if $\alpha_i\ne0$, \ $i=1,\dots,n$. By induction, we may prove this
for $i=j$ assuming that it is true for $i>j$. Then~\eqref{rpvgh} turns to the
identity
\begin{multline*}
  s_p\bigl( \v{G_1}^{\alpha_1}\cdots\v{G_j}^{\alpha_j} \bigr)
  (t_1\dots t_{j+1})^{\alpha_{j+1}}\cdots(t_1\cdots t_n)^{\alpha_n}\\
  =s_p\bigl(\v{H_1}^{\beta_1}\cdots\v{H_j}^{\beta_j}\bigr)
  (t_1\dots t_{j+1})^{\alpha_{j+1}}\cdots(t_1\cdots
  t_n)^{\alpha_n},
\end{multline*}
whence $s_p(\v{G_1}^{\alpha_1}\cdots\v{G_j}^{\alpha_j})=
s_p(\v{H_1}^{\beta_1}\cdots\v{H_j}^{\beta_j})$. Suppose that $\beta_l$ is the
last non-zero exponent (so that $\beta_{l+1}=\dots=\beta_j=0$). Then we also
have $\alpha_{l+1}=\dots=\alpha_j=0$, since otherwise
$s_p(\v{G_1}^{\alpha_1}\cdots\v{G_j}^{\alpha_j})$ would be divisible by
$t_1\dots t_{l+1}$, while $s_p(\v{H_1}^{\beta_1}\cdots\v{H_j}^{\beta_j})$ is
not. We also have $\alpha_l=\beta_l$ and $G_l=H_l$ since $\alpha_l$ is the
maximal power of $t_1\dots t_l$ that divides
$s_p(\v{G_1}^{\alpha_1}\cdots\v{G_j}^{\alpha_j})$. By induction, we conclude
that $\v{G_1}^{\alpha_1}\cdots\v{G_n}^{\alpha_n}=
\v{H_1}^{\beta_1}\cdots\v{H_n}^{\beta_n}$, whence $s_p(a)\ne0$.
\end{proof}

\begin{remark}
The same argument as in the proof of Lemma~\ref{alres} shows that
for arbitrary $Q$ the sum $s=\oplus_{{}_G}s_{{}_G}$ of (obviously
defined) restriction maps $s_{{}_G}$ over all minimal faces
$G\subset Q$ is a monomorphism.
\end{remark}

\begin{corollary}
The chain decomposition of $a\in\k[Q]$ is unique, and the monomials
$\v{G_1}^{\alpha_1}\cdots\v{G_q}^{\alpha_q}$ corresponding to all chains
$G_1\supset\cdots\supset G_q$ and all exponents $\alpha_i$ form an additive
basis of $\k[Q]$.
\end{corollary}

The \emph{$f$-vector} of $Q$ is defined as $\mb f(Q)=(f_0,\dots,f_{n-1})$
where $f_i$ is the number of faces of codimension $i+1$ (so that $f_0=m$ is
the number of facets).  The equivalent information is contained in the
\emph{$h$-vector} $\mb h(Q)=(h_0,\dots,h_n)$ determined by the equation
\begin{equation}
\label{hvector}
  h_0t^n+\ldots+h_{n-1}t+h_n=(t-1)^n+f_0(t-1)^{n-1}+\ldots+f_{n-1}.
\end{equation}
In particular, $h_0=1$ and
$h_n=(-1)^n+(-1)^{n-1}f_0+\dots+f_{n-1}$, which is equal to~$1$
when $Q$ is face-acyclic.

\begin{example}
We turn again to the $T^n$-action on $S^{2n}$ from Examples~\ref{2nsphere}
and~\ref{S2n orbit} and set $n=2$ there. Then $Q$ is a 2-ball with two
0-faces (say, $p$ and $q$) and two 1-faces (say, $G$ and $H$). Then $\mb
f(Q)=(2,2)$, $\mb h(Q)=(1,0,1)$ and
\[
  \k[Q]=\k[\v G,\v H,v_p,v_q]/(\v G\v H=v_p+v_q,\ v_pv_q=0),
\]
where $\deg\v G=\deg\v H=2$, $\deg v_p=\deg v_q=4$.
\end{example}

\subsection{Simplicial posets}\label{simplicial posets}

The faces (simplices) in a (finite) simplicial complex $K$ form a poset
(partially ordered set) with respect to the inclusion, and the empty simplex
$\varnothing$ is the initial element. This poset is called the \emph{face
poset} of~$K$, and it carries the same combinatorial information as the
simplicial complex itself. A poset $\mathcal P$ is called \emph{simplicial}
if it has an initial element $\hatzero$ and for each $x\in\mathcal P$ the
lower segment $[\hatzero,x]$ is a boolean lattice (the face poset of a
simplex). The face poset of a simplicial complex is a simplicial poset, but
there are simplicial posets that cannot be obtained in this way. In the
sequel we identify a simplicial complex with its face poset, thereby
regarding simplicial complexes as particular cases of simplicial posets.

To each $x\in\overline\P:=\mathcal P-\{\hatzero\}$ we assign a geometrical
simplex whose face poset is $[\hatzero,x]$, and glue these geometrical
simplices together according to the order relation in $\P$. We get a
cell complex such that the closure of each cell can be identified with a
simplex preserving the face structure and all the attaching maps are
inclusions. We call it a \emph{simplicial cell complex} and denote its
underlying space by $|\P|$. If $\P$ is (the face poset of) a simplicial
complex $K$, then $|\P|$ agrees with the geometric realisation $|K|$ of $K$.
The barycentric subdivision of a simplicial cell complex is obviously
defined, and is again a simplicial cell complex.

\begin{proposition}\label{baryc}
The barycentric subdivision of a simplicial cell complex is a (geometric
realisation of) simplicial complex.
\end{proposition}
\begin{proof}
Indeed, we may identify the barycentric subdivision under question with
the geometric realisation of the order complex $\Delta(\overline\P)$
of the poset $\overline\P$.
\end{proof}

In the sequel we will not distinguish between simplicial posets and
simplicial cell complexes, and call (the face poset of) the order complex
$\Delta(\overline\P)$ the \emph{barycentric subdivision} of $\P$.  The set of
faces of a nice manifold with corners $Q$ forms a simplicial poset with
respect to reversed inclusion (so $Q$ is the initial element).  We
call it the \emph{face poset} of $Q$.  It is a face poset of a simplicial
complex if and only if all non-empty multiple intersections of facets of $Q$
are connected.

\begin{example}
Let $Q$ be the orbit space from Example~\ref{S2n orbit}.  There are $n$
facets in $Q$ and the intersection of any $k$ facets is connected when $k\le
n-1$, but the intersection of $n$ facets consists of two points. The
corresponding simplicial cell complex is obtained by gluing two
$(n-1)$-simplices along their boundaries.
\end{example}

Let $\mathcal P$ be a simplicial poset.
When $[\hatzero,x]$ is the face poset of a $(k-1)$-simplex,
the rank of $x\in\overline{\mathcal P}$, denoted by
$\mathop{\mathrm{rk}}x=k$, is defined to be $k$.
The rank of $\mathcal P$ is the maximum of ranks of elements
in $\overline{\mathcal P}$.
Introduce the graded polynomial ring $\k[v_x\colon
x\in\overline{\mathcal P}]$ with $\deg v_x=2\mathop{\mathrm{rk}}x$.  We also
write formally $v_{\hatzero}=1$.  For any two elements $x,y\in\mathcal P$
denote by $x\vee y$ the set of their least common upper bounds, and by
$x\wedge y$ the set of their greatest common lower bounds. Since $\mathcal P$
is simplicial, $x\wedge y$ consists of a single element provided that $x\vee
y$ is non-empty.  The following is the obvious dualisation of
Definition~\ref{frgen}.

\begin{definition}
The \emph{face ring} of a simplicial poset $\mathcal P$ is the
quotient
\[
  \k[\mathcal P]:=
  \k[v_x\colon x\in\mathcal\overline{\mathcal P}]/\mathcal I_{\mathcal P},
\]
where $\mathcal I_{\mathcal P}$ is the ideal generated by the elements
\[
  v_xv_y-v_{x\wedge y}\cdot\!\!\!\sum_{z\in x\vee y}\!\!\! v_z.
\]
\end{definition}

\begin{remark}
Let $Q$ be a nice manifold with corners and let $\mathcal P$ be the face poset
of $Q$. Then $\k[Q]\cong \k[\mathcal P]$.  Let
$K$ be the nerve simplicial complex of the covering of
$\partial Q=\cup_{i=1}^m Q_i$ by the facets, that is, the simplicial
complex on $m$ vertices whose $(k-1)$-dimensional simplices correspond to the
codimension-$k$ prefaces of~$Q$.
If all non-empty multiple intersections of facets in $Q$ are connected,
then the Stanley--Reisner face ring $\k[K]$ agrees with $\k[\mathcal P]$,
but otherwise $\k[K]$ may differ from $\k[\mathcal P]$.
\end{remark}

The \emph{$f$-vector} of a simplicial poset $\mathcal P$ of rank $n$
is $\mb f(\mathcal
P)=(f_0,\dots,f_{n-1})$ where $f_i$ is the number of elements of rank~$i$.
The \emph{$h$-vector} $\mb h(\mathcal P)=(h_0,\dots,h_n)$ is determined
by~\eqref{hvector}. If $\P$ is the face poset of a nice manifold with corners
$Q$ then $\mb h(\P)=\mb h(Q)$.

Since we have defined $\deg v_x=2\mathop{\mathrm{rk}}x$, the face
ring $\k[\mathcal P]$ has no odd degree part.
Its Hilbert series $F(\k[\mathcal P];t):=\sum_i\dim_{\k}\k[\mathcal P]_{2i}
t^{2i}$, where $\k[\mathcal P]_{2i}$ denotes
the homogeneous degree $2i$ part of $\k[\mathcal P]$,
looks exactly as in the case of simplicial complexes.

\begin{theorem}[Proposition 3.8 of \cite{stan91}]
\label{psgfr}
Let $\mathcal P$ be a simplicial poset of rank $n$ with $h$-vector
$(h_0,h_1,\dots,h_n)$.  Then
$$
  F\bigl(\k[\mathcal P];t\bigr)
  =\frac{h_0+h_1t^2+\dots+h_nt^{2n}}{(1-t^2)^n}.
$$
\end{theorem}

In~\cite{da-ja91}, Davis and Januszkiewicz realised the classical
Stanley--Reisner face ring $\k[K]$ of a simplicial complex $K$ as the
equivariant cohomology ring of a $T$-space. The same approach works for a
simplicial poset $\P$ as well. The order complex $\Delta(\overline\P)$ is a
simplicial complex. Let $P$ be the cone on the geometric realisation
$|\Delta(\overline\P)|$. Since $|\Delta(\overline\P)|=|\P|$, the \lq\lq
boundary" of $P$ is $|\P|$.  For each simplex $\sigma\in
\Delta(\overline\P)$, let $F_\sigma\subset P$ denote the geometric
realisation of the poset $\{\tau\in\Delta(\overline\P)
\colon\sigma\subseteq\tau\}$. If $\sigma$ is a $(k-1)$-simplex, then we
declare $F_\sigma$ to be a \emph{face of codimension}~$k$. Therefore, each
facet (codimension-one face) can be identified with the star of some vertex
in $\Delta(\overline\P)$. Each codimension-$k$ face is a connected component
of an intersection of $k$ facets and is acyclic since it is a cone. In the
case when $\P$ is a simplicial complex the space $P$ with the face
decomposition was called in~\cite[p.~428]{da-ja91} a \emph{simple polyhedral
complex}.

Suppose that the number of facets of $P$ is $m$ and that we have a map $\L$
as in (\ref{eqn:characteristic map}) satisfying the condition form
Lemma~\ref{nonsi}. (The existence of such a map $\L$ is equivalent to the
existence of a \emph{linear system of parameters} in the ring $\Z[\mathcal
P]$, see e.g.~\cite[Lemma~III.2.4]{stan96}.) Then the same construction as
$M_Q(\L)$ in (\ref{idspa}) with $Q$ replaced by $P$ produces a $T$-space
$M_P(\L)$.  Since $P$ is not a manifold with corners for arbitrary~$\P$, the
space $M_P(\L)$ may fail to be a manifold.  Nevertheless, a similar argument
to that in~\cite[Theorem 4.8]{da-ja91} gives the following result:

\begin{proposition} \label{prop:srfr}
$H^*_T(M_P(\L);\Z)$ is isomorphic to $\Z[\P]$ as a ring.
\end{proposition}

For an arbitrary nice manifold with corners $Q$ the equivariant cohomology of
the canonical model $M_Q(\L)$ may fail to be isomorphic to $\Z[Q]$ as the
faces of $Q$ themselves may have complicated cohomology. In the next sections
we shall study this question in more details. As the first step in this
direction we relate $M_Q(\L)$ to $M_P(\L)$ in our last statement of this
paragraph.

\begin{proposition}\label{eqn:Phi}
Let $Q$ be a nice manifold with corners, and $P$ the space associated with
the face poset $\P$ of~$Q$. Then there is a map $Q\to P$ which preserves the
face structure. It is covered by a canonical equivariant map
\[
  \Phi \colon M_Q(\L)\to M_P(\L).
\]
\end{proposition}
\begin{proof}
The map $Q\to P$ is constructed inductively, starting from an identification
of vertices and extending the map on each higher-dimensional face by a
degree-one map. Every face of $P$ is a cone, so there are no obstructions to
such extensions. Since the map between orbit spaces preserves the face
structure, it is covered by an equivariant map of the identification spaces
\[
  M_Q(\L)=T\times Q/\!\sim\;\longrightarrow T\times P/\!\sim\:=M_P(\L)
\]
by the definition of identification spaces, see~\eqref{idspa}.
\end{proof}

\section{Axial functions and Thom classes}

Here we relate the equivariant cohomology ring of a torus manifolds $M$ to
the face ring of the orbit space~$Q$.  We construct a natural ring
homomorphism from $\Z[Q]$ to $H^*_T(M)$ modulo $H^*(BT)$-torsions. In the
next section we show that this is an isomorphism when $H^{odd}(M)=0$.
In this and next sections we assume that $M$ is locally standard for
simplicity, but the arguments will work without this assumption with a
little modification.

\subsection{Axial functions}

Like in the algebraic situation of the previous section, we have the
restriction map to a sum of polynomial rings:
\begin{equation}\label{eqn:restriction map}
  r=\bigoplus_{p\in M^T} r_p\colon H^*_T(M)\to H^*_T(M^T)=
  \bigoplus_{p\in M^T}H^*(BT).
\end{equation}
The kernel of $r$ is the $H^*(BT)$-torsion subgroup of $H^*_T(M)$, so
$r$ is injective when $H^{odd}(M)=0$ by Lemma~\ref{lemm:freeness}.

We identify $M^T$ with the vertices of $Q$. The $1$-skeleton of $Q$,
consisting of vertices ($0$-faces) and edges ($1$-faces) of $Q$, is an
$n$-valent graph. Denote by $E(Q)$ the set of oriented edges.  Given an
element $e\in E(Q)$, denote the initial point and terminal point of $e$ by
$i(e)$ and $t(e)$ respectively. Then $\Me:=\pi^{-1}(e)$ is a 2-sphere fixed
by a codimension-one subtorus in $T$ (here $\pi\colon M\to Q$ is the quotient
map). It contains two $T$-fixed points $i(e)$ and $t(e)$. The 2-dimensional
subspace $\mathcal T_{i(e)}\Me\subseteq \mathcal T_{i(e)}M$ is an irreducible
component of the tangential $T$-representation $\mathcal T_{i(e)}M$. The same
is true for the other $T$-fixed point $t(e)$, and the $T$-representations
$\mathcal T_{i(e)}M$ and $\mathcal T_{t(e)}M$ are isomorphic. There is a
unique characteristic submanifold, say $M_i$, intersecting $\Me$ at $i(e)$
transversally.  Assuming both $M$ and $M_i$ are oriented, we choose a
compatible orientation for the normal bundle $\nu_i$ of $M_i$ and therefore,
for $\mathcal T_{i(e)}\Me$. The orientation on $\mathcal T_{i(e)}\Me$
determines a complex structure, so that $\mathcal T_{i(e)}\Me$ can be viewed
as a complex $1$-dimensional $T$-representation.  This defines an element of
$\Hom(T,S^1)=H^2(BT)$, which we denote by $\alpha(e)$.

Let $e^T(\nu_i)$ be the equivariant Euler class in $H^2_T(M_i)$ and denote
its restriction to $p\in M_i^T$ by $e^T(\nu_i)|_p\in H^2_T(p)=H^2(BT)$. Then
\begin{equation} \label{eqn:euler}
  e^T(\nu_i)|_p=\alpha(e),
\end{equation}
where $e$ is the unique edge such that $i(e)=p$ and $e\notin \F_i=\pi(M_i)$.
Following~\cite{GZ}, we call the map
\[
  \alpha\colon E(Q)\to H^2(BT)
\]
an \emph{axial function}.

\begin{lemma}
The axial function $\alpha$ has the following properties:
\begin{itemize}
\item[(1)] $\alpha(\bar e)=\pm \alpha(e)$ for all $e\in E(Q)$, where $\bar e$
denotes $e$ with the opposite orientation;
\item[(2)] for each vertex (or a $T$-fixed point) $p$, the set $\alpha_p:=\{
\alpha(e)\colon i(e)=p\}$ is a basis of $H^2(BT)$ over $\Z$.
\item[(3)] for $e\in E(Q)$, we have
$\alpha_{i(e)}\equiv \alpha_{t(e)}\mod \alpha(e)$.
\end{itemize}
\end{lemma}
\begin{proof}
Property (1) follows from the fact that $\mathcal T_{i(e)}\Me$ and $\mathcal
T_{t(e)}\Me$ are isomorphic as real $T$-representations, and~(2) from that
the $T$-representation $\mathcal T_{i(e)}M$ is faithful of complex
dimension~$n$. Let $T_e$ be the codimension one subtorus fixing $\Me$. Then
the $T$-representations $\mathcal T_{i(e)}M$ and $\mathcal T_{t(e)}M$ are
isomorphic as $T_e$-representations, since the points $i(e)$ and $t(e)$ are
contained in the same connected component $\Me$ of the $T_e$-fixed point set.
This implies~(3).
\end{proof}

\begin{remark}
In \cite{GZ}, the property $\alpha(\bar e)=-\alpha(e)$ is required in the
definition of axial function, but we allow $\alpha(\bar e)=\alpha(e)$. For
example, $\alpha(\bar e)=\alpha(e)$ for the $T^2$-action on $S^{4}$ from
Example~\ref{2nsphere}.
\end{remark}

\begin{lemma} \label{lemm:necessity}
Fix $\eta\in H^*_T(M)$; then $r_{i(e)}(\eta)-r_{t(e)}(\eta)$ is divisible by
$\alpha(e)$ for all $e\in E(Q)$.
\end{lemma}
\begin{proof}
Consider the commutative diagram of restrictions
\[
\begin{CD}
  H^*_T(M) @>>> H^*_T(i(e))\oplus H^*_T(t(e))= @. H^*(BT)\oplus H^*(BT)\\
  @VVV          @VVV \\
  H^*_{T_e}(\Me) @>>> H^*_{T_e}(i(e))\oplus H^*_{T_e}(t(e))= @.\:
  H^*(BT_e)\oplus H^*(BT_e)
\end{CD}.
\]
Since $H^*_{T_e}(\Me)=H^*(BT_e)\otimes H^*(\Me)$, the two components of the
image of $\eta$ in $H^*(BT_e)\oplus H^*(BT_e)$ above coincide. Therefore it
follows from the commutativity of the above diagram that the restrictions of
$r_{i(e)}(\eta)$ and $r_{t(e)}(\eta)$ to $H^*(BT_e)$ coincide. Since the
kernel of the restriction map $H^*(BT)\to H^*(BT_e)$ is the ideal generated
by $\alpha(e)$, the lemma follows.
\end{proof}

\subsection{Thom classes}

The preimage $\MG:=\pi^{-1}(\G)$ of a codimension-$k$ face $\G\subset Q$ is a
connected component of an intersection of $k$ characteristic submanifolds.
The orientations of $M$ and characteristic submanifolds $M_i$ determine
compatible orientations for the normal bundles $\nu_i$ of $M_i$. These
orientations determine an orientation on the normal bundle $\n{\G}$ of $\MG$,
and thereby on $\MG$ itself, since $M$ is oriented. With this convention on
orientations, we consider the Gysin homomorphism $H_T^0(\MG)\to H_T^{2k}(M)$
in the equivariant cohomology and denote the image of the identity element by
$\ta{\G}$. The element $\ta{\G}$ may be thought of as the Poincar\'e dual of
$\MG$ in equivariant cohomology and is called the \emph{Thom class} of $\MG$.
The restriction of $\ta{\G}\in H_T^{2k}(M)$ to $H_T^{2k}(\MG)$ is the
equivariant Euler class of $\n{\G}$, and $r_p(\ta{\G})=0$ unless
$p\in(\MG)^T$. It follows from (\ref{eqn:euler}) that
\begin{equation} \label{eqn:tau_G}
  r_p(\ta{\G})=\left\{%
  \begin{array}{ll}\displaystyle
    \prod_{i(e)=p,\ e\nsubseteq \G}\alpha(e),
       & \hbox{if \ $p\in (\MG)^T$;}\\
   \qquad 0, & \hbox{otherwise.}\\
  \end{array}%
  \right.
\end{equation}

We set
\[
  \H^*_T(M):=H^*_T(M)/H^*(BT)\text{-torsions}.
\]
The restriction map~(\ref{eqn:restriction map}) induces a monomorphism
$\H^*_T(M)\to H^*_T(M^T)$, which we also denote by $r$.
Therefore, $\ta{\G}=0$ in $\H^*_T(M)$ if $M_F$ has no $T$-fixed point.
The following lemma
shows that the relations from Definition~\ref{frgen} hold in $\H^*_T(M)$ with
$\v{\G}$ replaced by $\ta{\G}$.

\begin{lemma}
For any two faces $G$ and $H$ of $Q$, the relation
\[
  \ta{G}\ta{H}=
    \ta{G\vee H}\cdot\!\!\!\sum_{E\in{G\cap H}}\!\!\!\ta{E},
\]
holds in $\H^*_T(M)$, where we set $\tau_{\varnothing}=0$.
\end{lemma}
\begin{proof}
Since the restriction map $r\colon \H^*_T(M)\to H^*_T(M^T)$ is injective, it
suffices to show that $r_p$ maps both sides of the identity to the same
element for all $p\in M^T$.

Let $p\in M^T$. For a face $\G$ such that $p\in\G$, we set
\[
  N_p(\G):= \{ e\in E(Q) \colon i(e)=p,\ e\notin \G\},
\]
which may be thought of as the set of directions normal to $\G$ at $p$.  We
also set $N_p(\G)=\varnothing$ if $p\notin\G$.  Then the identity
(\ref{eqn:tau_G}) can be written as
\begin{equation} \label{eqn:tau_G2}
  r_p(\ta{\G})=\prod_{e\in N_p(\G)}\alpha(e)
\end{equation}
where the right hand side is understood to be zero if $N_p(\G)=\varnothing$.
If $p\notin G\cap H$, then $N_p(E)=\varnothing$ for any connected component
$E$ of $G\cap H$ and either $N_p(G)=\varnothing$ or $N_p(H)=\varnothing$.
Therefore, both sides of the identity from the lemma map to zero by~$r_p$.
If $p\in G\cap H$, then
\[
  N_p(G)\cup N_p(H)=N_p(G\vee H)\cup N_p(E)
\]
where $E$ is the connected component of $G\cap H$ containing $p$, and
$N_p(E')=\varnothing$ for any other connected component of $G\cap H$. This
together with (\ref{eqn:tau_G2}) shows that both sides of the identity
map to the same element by~$r_p$.
\end{proof}

By virtue of the above lemma, the map $\Z[\v{\G}\colon \G \text{ a face}]\to
H^*_T(M)$ sending $\v{\G}$ to $\ta{\G}$ induces a homomorphism
\begin{equation} \label{eqn:varphi2}
  \f\colon \Z[Q]\to \H_T^*(M).
\end{equation}

\begin{lemma} \label{lemm:monof}
The homomorphism $\f$ is injective if every face of $Q$ has a vertex.
\end{lemma}

\begin{proof}
We have $s=r\circ \f$, where $s$ is the map from Lemma~\ref{alres}.
Since $s$ is injective if every face of $Q$ has a vertex, so is~$\f$.
\end{proof}

\section{Equivariant cohomology ring of torus manifolds with vanishing
odd-degree cohomology}

In this section we give a sufficient condition for the monomorphism $\f$ in
(\ref{eqn:varphi2}) to be an isomorphism (Theorem~\ref{theo:eqcoh}).  In
particular, it turns out that $\f$ is an isomorphism when $H^{odd}(M)=0$
(Corollary~\ref{coro:eqcoh}). Using these results, we give a description of
the ring structure in $H^*(M)$ in the case when $H^{odd}(M)=0$
(Corollary~\ref{coro:stcoh}).

\subsection{Ring structure in equivariant cohomology}\label{rings}
The following theorem shows that the converse of Lemma~\ref{lemm:necessity}
holds for torus manifolds with vanishing odd degree cohomology.

\begin{theorem}[\cite{g-k-m98}, see also Chapter 11 in \cite{GS}]
\label{theo:GKM}
Suppose $H^{odd}(M)=0$ and we are given an element $\eta_p\in H^*(BT)$ for
each $p\in M^T$.  Then $(\eta_p)\in \bigoplus_{p\in M^T}H^*(BT)$ belongs to
the image of the restriction map $r$ in~\eqref{eqn:restriction map} if and
only if $\eta_{i(e)}-\eta_{t(e)}$ is divisible by $\alpha(e)$ for any $e\in
E(Q)$.
\end{theorem}

\begin{corollary}\label{1skel}
The $1$-skeleton of any face of $Q$ (including $Q$ itself) is connected if
$H^{odd}(M)=0$.
\end{corollary}
\begin{proof}
Since $M$ is connected, the image $r(H^0_T(M))$ is one-dimensional. Then it
follows from the ``if" part of Theorem~\ref{theo:GKM} that the 1-skeleton of
$Q$ is connected. Similarly, the 1-skeleton of any face $\G$ of $Q$ is
connected
because $\MG=\pi^{-1}(\G)$ is also a torus manifold with vanishing odd degree
cohomology (see Lemma~\ref{lemm:odd=0}).
\end{proof}

\begin{remark}
The connectedness of $1$-skeletons of faces of $Q$ can be proven without
referring to Theorem~\ref{theo:GKM}, see remark after Theorem~\ref{cohfa}.
\end{remark}

For a face $\G$ of $Q$, we denote by $I(\G)$ the ideal in $H^*(BT)$ generated
by all elements $\alpha(e)$ with $e\in \G$.

\begin{lemma} \label{lemm:I(G)}
Suppose that the $1$-skeleton of a face $\G$ is connected. Given $\eta\in
H^*_T(M)$, if $r_p(\eta)\notin I(\G)$ for some vertex $p\in \G$, then
$r_q(\eta)\notin I(\G)$ for any vertex $q\in \G$.
\end{lemma}
\begin{proof}
Suppose $r_q(\eta)\in I(\G)$ for some vertex $q\in \G$. Then $r_s(\eta)\in
I(\G)$ for any vertex $s\in \G$ joined to $q$ by an edge $f\subseteq\G$
because $r_q(\eta)-r_s(\eta)$ is divisible by $\alpha(f)$ by
Lemma~\ref{lemm:necessity}. Since the $1$-skeleton of $\G$ is connected,
$\eta(q)\in I(\G)$ for any vertex $q\in \G$, which contradicts the
assumption.
\end{proof}

\begin{proposition}\label{prop:modgn}
If the $1$-skeleton of
every face of $Q$ is connected, then $\H^*_T(M)$ is generated by the elements
$\ta{\G}$ as an $H^*(BT)$-module.
\end{proposition}
\begin{proof}
Let $\eta\in H^{>0}_T(M)$ be a nonzero element. Set
$$
  Z(\eta):=\{p\in M^T\colon r_p(\eta)=0\}.
$$
Take $p\in M^T$ such that $p\notin Z(\eta)$. Then $r_p(\eta)\in H^*(BT)$ is
non-zero and we can express it as a polynomial in $\{\alpha(e)\colon
i(e)=p\}$ (the latter is a basis of $H^2(BT)$).  Let
\begin{equation}
\label{monom}
  \prod_{i(e)=p}\alpha(e)^{n_e},
\end{equation}
$n_e\ge0$, be a monomial entering $r_p(\eta)$ with a non-zero coefficient.
Let $\G$ be the face spanned by the edges $e$ with $n_e=0$. Then
$r_p(\eta)\notin I(\G)$ since $r_p(\eta)$ contains the
monomial~\eqref{monom}. Hence, $r_q(\eta)\notin I(\G)$, in particular
$r_q(\eta)\ne0$, for every vertex $q\in \G$ by Lemma~\ref{lemm:I(G)}.

On the other hand, it follows from (\ref{eqn:tau_G}) that the monomial
\eqref{monom} can be written as $r_p(\u{\G}\ta{\G})$ with some $\u{\G}\in
H^*(BT)$. Set $\eta':=\eta-\u{\G}\ta{\G}\in H^*_T(M)$.  Since
$r_q(\ta{\G})=0$ for every vertex $q\notin \G$, we have
$r_q(\eta')=r_q(\eta)$ for such $q$.  At the same time, $r_q(\eta)\not=0$ for
every vertex $q\in \G$ (see above).  It follows that $Z(\eta')\supseteq
Z(\eta)$. However, the number of monomials in $r_p(\eta')$ is less than that
in $r_p(\eta)$. Therefore, subtracting from $\eta$ a linear combination of
$\ta{\G}$'s with coefficients in $H^*(BT)$, we obtain an element $\lambda$
such that $Z(\lambda)$ contains $Z(\eta)$ as a proper subset. Repeating this
procedure, we end up at an element whose restriction to every vertex is zero.
Since the restriction map $r\colon \H^*_T(M) \to H^*_T(M^T)$ is injective,
this finishes the proof.
\end{proof}

\begin{theorem}\label{theo:eqcoh}
Let $M$ be a (locally standard) torus manifold with orbit space $Q$.
If every face of $Q$ has a vertex and its $1$-skeleton is connected, then the
mo\-no\-mor\-phism $\f\colon \Z[Q]\to \H^*_T(M)$ in \eqref{eqn:varphi2} is an
isomorphism.
\end{theorem}
\begin{proof}
To prove that $\f$ is surjective it suffices to show that $\H^*_T(M)$ is
generated by the elements $\ta{\G}$ as a ring. By Proposition~\ref{prop:a_i},
$\H^2_T(M)$ is generated over $\Z$ by the elements $\ta{\F_i}$ corresponding
to the facets $\F_i$. (Note: the notation $\tau_i$ is used for $\ta{\F_i}$ in
Proposition~\ref{prop:a_i}.) In particular, any element in $H^2(BT)\subset
\H^*_T(M)$ can be written as a linear combination of $\ta{\F_i}$'s with
coefficients in~$\Z$.  Hence, any element in $H^*(BT)$ is a polynomial in
$\ta{\F_i}$'s. The rest follows from Proposition~\ref{prop:modgn}.
\end{proof}

Now assume $H^{odd}(M)=0$. Then the assumption in Theorem~\ref{theo:eqcoh}
is satisfied, and $H^*_T(M)$ is a free $H^*(BT)$-module by
Lemma \ref{lemm:freeness}, whence $\H^*_T(M)=H^*_T(M)$.

\begin{corollary}\label{coro:eqcoh}
For a torus manifold $M$ with vanishing odd degree cohomology, the map
$\f\colon \Z[Q]\to H^*_T(M)$ in \eqref{eqn:varphi2} is an isomorphism.
\end{corollary}
\begin{proof}
This follows from Corollary~\ref{1skel} and Theorem~\ref{theo:eqcoh}.
\end{proof}

\begin{remark}
When $H^*(M)$ is generated in degree two, all non-empty multiple
intersections of facets are connected by Lemma~\ref{cmcon}. Therefore, the
face poset of $Q$ is the face poset of the nerve $K$ of the covering of
$\partial Q$, and $\Z[Q]$ reduces to the classical Stanley--Reisner face ring
of a simplicial complex. Therefore, Corollary~\ref{coro:eqcoh} is a
generalisation of Proposition 3.4 in \cite{masu99}.
\end{remark}

If $\P$ is the face poset of $Q$, then $\Z[\P]=\Z[Q]$ by the definition.
The following statement gives a characterisation of torus manifolds $M$ with
vanishing odd degree cohomology (and with cohomology generated in degree two)
in terms of the face poset $\P$ associated with $M$.

\begin{theorem} \label{theo:face CM}
Let $M$ be a torus manifold with orbit space $Q$, and let $\P$ be the face
poset of $Q$.  Then $H^{odd}(M)=0$ if and only if the following two
conditions are satisfied:
\begin{itemize}
\item[(1)] $H^*_T(M)$ is isomorphic to $\Z[\P](=\Z[Q])$ as a ring;
\item[(2)] $\Z[\P]$ is Cohen--Macaulay.
\end{itemize}
Moreover, $H^*(M)$ is generated by its degree-two part if and only if $\P$ is
(the face poset of) a simplicial complex in addition to the above two
conditions.
\end{theorem}
\begin{proof}
If $H^{odd}(M)=0$, then $H^*_T(M)\cong \Z[Q]=\Z[\P]$ by
Corollary~\ref{coro:eqcoh}, and $\Z[\P]$ is Cohen--Macaulay because
$H^*_T(M)$ is a free $H^*(BT)$-module by Lemma~\ref{lemm:freeness}.  This
proves the \lq\lq only if" part of the first statement.

Now we prove the ``if" part. Let $\l\colon ET\times_T M\to BT$ be the
projection, and consider the composite map
\[
  H^*(BT)\stackrel{\l^*}\longrightarrow  H^*_T(M)
  \stackrel{r}\longrightarrow \bigoplus_{p\in M^T}H^*(BT).
\]
Its restriction to each summand of the target is the identity, i.e., $r\circ
\l^*$ is a diagonal map. This implies that $\l^*(t_1),\dots,\l^*(t_n)$ is a
linear system of parameters (an l.s.o.p.), see
\cite[Theorem~5.1.16]{br-he98}. By the assumption, $H^*_T(M)$ is isomorphic
to $\Z[\P]$ and $\Z[\P]$ is Cohen--Macaulay, so every l.s.o.p. is a regular
sequence (see \cite[Theorem I.5.9]{stan96}). It follows that $H^*_T(M)$ is a
free $H^*(BT)$-module and hence $H^{odd}(M)=0$ by Lemma~\ref{lemm:freeness},
thus proving the ``if" part of the first statement.

It remains to prove the second statement. The \lq\lq only if" part follows
from Lemma~\ref{cmcon} by the last remark. For the ``if" part, if
$\P$ is a simplicial poset, then $\Z[\P]$ is generated by
its degree-two part. By the first statement of the theorem, $H^*_T(M)\cong
\Z[\P]$ is a free $H^*(BT)$-module, whence $H^*(M)$ is a quotient ring of
$H^*_T(M)$. It follows that $H^*(M)$ is also generated by its degree-two
part.
\end{proof}

The following description of cohomology ring of a torus manifold with
vanishing odd degree cohomology generalises that of a complete
non-singular toric variety, see \cite[p.106]{fult93}.

\begin{corollary} \label{coro:stcoh}
For a torus manifold $M$ with vanishing odd degree cohomology,
\[
  H^*(M)\cong \Z[\v{\G}\colon \G
  \text{ a face of $Q$}]/I \qquad\text{as a ring,}
\]
where $I$ is the ideal generated by the following two types of elements:
\begin{itemize}
\item[(1)]  $\displaystyle{\v{G}\v{H}-
    \v{G\vee H}\sum_{E\in{G\cap H}}\v{E}} ;
$
\item[(2)] $\displaystyle{\sum_{i=1}^m\< t,a_i\> \v{\F_i}}$
for $t\in H^2(BT)$.
\end{itemize}
Here $\F_i$ are the facets of $Q$ and the elements $a_i\in H_2(BT)$ are
defined in Proposition~\ref{prop:a_i}.
\end{corollary}
\begin{proof}
Since the Serre spectral sequence of the fibration $\l\colon ET\times_T M\to
BT$ collapses, the restriction map $H^*_T(M)\to H^*(M)$ is surjective and its
kernel is the ideal generated by all $\l^*(t)$ with $t\in H^2(BT)$.
Therefore, the statement follows from Proposition~\ref{prop:a_i} and
Corollary~\ref{coro:eqcoh}.
\end{proof}

\subsection{Dehn-Sommerville equations}\label{sub:DS}

Suppose that $H^{odd}(M)=0$. Then, since $H^*_T(M)=H^*(BT)\otimes H^*(M)$ by
Lemma~\ref{lemm:freeness} and $H^*(BT)$ is a polynomial ring in $n$ variables
of degree two, the Hilbert series of $H^*_T(M)$ is given by
\[
  F\bigl(H^*_T(M);t\bigr)=\frac{\sum_{i=0}^n \rank_\Z
  H^{2i}(M)t^{2i}}{(1-t^2)^n}.
\]
On the other hand, the Hilbert series of the face ring $\Z[Q]$ is given by
Theorem~\ref{psgfr} and these two series must coincide by
Corollary~\ref{coro:eqcoh}. It follows that
\begin{equation}\label{betti}
  \rank_\Z H^{2i}(M)=h_i.
\end{equation}
Since $M$ is a manifold, the Poincar\'e duality implies that
\begin{equation}\label{dsequ}
  h_i=h_{n-i},\quad i=0,\dots,n.
\end{equation}
When every non-empty multiple intersection of facets in $Q$ is connected,
$\Z[Q]$ reduces to the classical Stanley--Reisner ring of the nerve of the
covering of $\partial Q$ and equations (\ref{dsequ}) are known as the
\emph{Dehn--Sommerville equations} for the numbers of faces.

\section{Orbit spaces of torus manifolds with cohomology generated in degree
two}

Using the equivariant cohomology calculations from the previous section, we
are finally able to relate the cohomology of a torus manifold $M$ and the
cohomology of its orbit space $Q$. The main result of this section is
Theorem~\ref{hptpe} which gives a cohomological characterisation of torus
manifolds whose orbit spaces are homology polytopes. Using this result, in
the next section we prove that $Q$ is face-acyclic if $H^{odd}(M)=0$.

\begin{lemma}\label{lemm:H1=0}
If $H^{odd}(M)=0$, then $H^1(Q;\k)=0$ for any coefficient ring $\k$.
In particular, $Q$ is orientable.
\end{lemma}
\begin{proof}
We use the Leray spectral sequence (with $\k$ coefficient) of the projection
map $ET\times_T M\to M/T=Q$ on the second factor.  Its $E_2$ term is given by
$E_2^{p,q}=H^p(M/T;\mathcal H^q)$ where $\mathcal H^q$ is a sheaf with stalk
$H^q(BT_x;\k)$ over a point $x\in M/T$, and the spectral sequence converges
to $H^*_T(M;\k)$. Since the $T$-action on $M$ is locally standard by
Theorem~\ref{theo:local standardness}, the isotropy group $T_x$ at $x\in M$
is a subtorus; so $H^{odd}(BT_x;\k)=0$. Hence, $\mathcal
H^{odd}=0$, in particular, $\mathcal H^1=0$. Moreover, $\mathcal H^0=\k$ (a
constant sheaf). Therefore, we have $E_2^{0,1}=0$ and
$E_2^{1,0}=H^1(M/T;\k)$, whence $H^1(M/T;\k)\cong H^1_T(M;\k)$.  On the other
hand, since $H^{odd}(M)=0$ by assumption, $H^*_T(M)$ is a free
$H^*(BT)$-module (isomorphic to $H^*(BT)\otimes H^*(M)$ by
Lemma~\ref{lemm:freeness}). Therefore, $H^{odd}_T(M;\k)=0$ by the universal
coefficient theorem. In particular, $H^1_T(M;\k)=0$, thus proving the lemma.
\end{proof}

\begin{lemma}\label{kgore}
If either
\begin{itemize}
\item[(1)] $Q$ is a homology polytope, or
\item[(2)] $H^*(M)$ is generated by its degree-two part,
\end{itemize}
then the face poset $\P$ of $Q$ is (the face poset of) a simplicial
Gorenstein* complex. In particular, $\Z[\P]$ is Cohen--Macaulay and the
geometric realisation $|\P|$ of $\P$ has the homology of an $(n-1)$-sphere.
\end{lemma}
\begin{proof}
Under either assumption (1) or (2), all non-empty multiple intersections of
facets of $Q$ are connected, so $\P$ agrees with the face poset of the nerve
simplicial complex $K$ of the covering of $\partial Q$.  In what follows we
identify $\P$ with $K$.

First we prove that $\P$ is Gorenstein* under assumption~(1).  According to
Theorem II.5.1 of \cite{stan96} it is enough to show that the link of a
simplex $\sigma$ of $\K$, denoted by $\link\sigma$, has the homology of a
sphere of $\dim \link\sigma=n-2-\dim\sigma$.  If $\sigma=\varnothing$ then
$\link\sigma$ is $\K$ itself and its homology is isomorphic to the homology
of the boundary $\partial Q$ of $Q$, since $\K$ is the nerve of $Q$ and $Q$
is a homology polytope. If $\sigma\not=\varnothing$ then $\link\sigma$ is the
nerve of a face of $Q$.  Since any face of $Q$ is again a homology polytope,
$\link\sigma$ has the homology of a sphere of $\dim\link\sigma$ by the same
argument.

Now we prove that $\K$ is Gorenstein* under assumption~(2). Using
Theorem~II.5.1 of~\cite{stan96} once again, it is enough to show that
\begin{itemize}
\item[(a)] $\K$ is Cohen--Macaulay;
\item[(b)] every $(n-2)$-dimensional simplex is contained in exactly two
$(n-1)$-\-di\-men\-si\-o\-nal simplices;
\item[(c)] $\chi(\K)=\chi(S^{n-1})$.
\end{itemize}
The condition (a) follows from Lemma~\ref{lemm:freeness} and
Corollary~\ref{coro:eqcoh}. By definition, every $k$-dimensional simplex of
$\K$ corresponds to a set of $k+1$ characteristic submanifolds having
non-empty intersection. By Lemma~\ref{cmcon}, the intersection of any $n$
characteristic submanifolds is either empty or consists of a single $T$-fixed
point. This means that the $(n-1)$-simplices of $\K$ are in one-to-one
correspondence with the $T$-fixed points of $M$.  Now, each $(n-2)$-simplex
of $\K$ corresponds to a non-empty intersection of $n-1$ characteristic
submanifolds of~$M$. The latter intersection is connected by
Lemma~\ref{cmcon} and has a non-trivial $T$-action, so it is a $2$-sphere.
Every $2$-sphere contains exactly two $T$-fixed points, which implies~(b).
Finally,~(c) is just the Dehn--Sommerville equation $h_0=h_n$,
see~\eqref{hvector} and~\eqref{dsequ}.
\end{proof}

\begin{theorem}\label{hptpe}
The cohomology of a torus manifold $M$ is generated by its degree-two part if
and only if $M$ is locally standard and the orbit space $Q$ is a homology
polytope.
\end{theorem}
\begin{proof}
Let $\P$ be the face poset of $Q$, and $P$ the cone on~$|\P|$
with the face structure associated with $\P$, see end of
subsection~\ref{simplicial posets}.

We first prove the ``if" part. Suppose $Q$ is a homology polytope. Since
$H^2(Q)=0$ and $M$ is locally standard,
$M$ is equivariantly homeomorphic to $M_Q(\L)$ by
Lemma~\ref{acydj}; so we may regard the map $\Phi$ in (\ref{eqn:Phi}) as a
map from $M$ to $M_P:=M_P(\L)$. Let $M_{P,i}$ be characteristic subcomplexes
of $M_P$ defined similarly to characteristic submanifolds $M_i$ of $M$.
Since the $T$-actions on $M_P\backslash \cup_i M_{P,i}$ and $M\backslash
\cup_i M_{i}$ are free, we have
\begin{equation*} \label{eqn:relative}
  H_T^*(M_P,\cup_i M_{P,i})\cong H^*(P,|\P|), \quad
  H^*_T(M,\cup_i M_{i})\cong H^*(Q,\partial Q).
\end{equation*}
Therefore, the map $\Phi$ induces a map between exact sequences
\begin{equation}\label{seqma}
\begin{CD}
  \longrightarrow\; H^*(P,|\P|) @>>> H_T^*(M_P) @>>>
  H_T^*(\cup_i M_{P,i})\;\longrightarrow\\
  @VVV @VV\Phi^*V @VVV\\
  \longrightarrow\; H^*(Q,\partial Q) @>>> H_T^*(M) @>>>
  H_T^*(\cup_i M_{i})\;\longrightarrow
\end{CD}
\end{equation}
Each $M_{i}$ itself is a torus manifold over a homology polytope $\F_i$.
Using induction and a Mayer--Vietoris argument, we may assume that the map
$H_T^*(\cup_i M_{P,i})\to H_T^*(\cup_i M_{i})$ above is an isomorphism.  By
Lemma~\ref{kgore}, $|\P|$ has the homology of an $(n-1)$-sphere, and since
$P$ is the cone over $|\mathcal P|$, we have $H^*(P,|\P|)\cong
H^*(D^n,S^{n-1})$.  We also have $H^*(Q,\partial Q)\cong H^*(D^n,S^{n-1})$
because $Q$ is a homology polytope. Using these isomorphisms, we see from the
construction of the map $\Phi$ that the induced map $H^*(P,|\P|)\to
H^*(Q,\partial Q)$ is the identity map on $H^*(D^n,S^{n-1})$. Therefore, the
5-lemma applied to~\eqref{seqma} shows that $\Phi^*\colon H_T^*(M_P)\to
H_T^*(M)$ is an isomorphism; whence $H_T^*(M)\cong \Z[\P]$ by
Proposition~\ref{prop:srfr}.  We also know that $\Z[\K]$ is Cohen--Macaulay
by Lemma~\ref{kgore}.  Therefore, the two conditions in
Theorem~\ref{theo:face CM} are satisfied.  It follows that $H^*(M)$ is
generated by its degree-two part by Theorem~\ref{theo:face CM}, which
finishes the proof of the ``if" part.

\smallskip

Now we prove the ``only if" part. Suppose that $H^*(M)$ is generated by the
degree-two elements.  Then $M$ is locally standard by Theorem~\ref{theo:local
standardness}.  Since all non-empty multiple intersections of characteristic
submanifolds are connected and their cohomology rings are generated in degree
two by Lemma~\ref{cmcon}, we may assume by induction that all the proper
faces of $Q$ are homology polytopes. In particular, the proper faces are
acyclic, whence $H^*(\partial Q)\cong H^*(|\P|)$.  This together with
Lemma~\ref{kgore} shows that
\begin{equation}\label{eqn:boundary Q}
  H^*(\partial Q)\cong H^*(S^{n-1}).
\end{equation}

\smallskip
\noindent
{\bf Claim.}
$H^2(Q)=0$.

\smallskip
\noindent
The claim is trivial for $n=1$. If $n=2$ then $Q$ is a surface with boundary,
hence, $H^2(Q)=0$ in this case too. Now assume $n\ge 3$.  Let us consider the
exact equivariant cohomology sequence of pair $(M,\cup_i M_i)$, see
the bottom row of~\eqref{seqma}. All the maps in the exact sequence are
$H^*(BT)$-module maps. By Lemma~\ref{lemm:freeness}, $H^*_T(M)$ is a free
$H^*(BT)$-module. On the other hand, $H^*(Q,\partial Q)$ is finitely
generated over $\Z$, so it is a torsion $H^*(BT)$-module.  It follows that
the whole sequence splits in short exact sequences:
\begin{equation}\label{short exact}
  0 \to H^k_T(M) \to H^k_T(\cup_i M_i) \to H^{k+1}(Q,\partial Q) \to 0
\end{equation}
Taking $k=1$ above, we get
\[
  H^1_T(\cup_iM_i)\cong H^2(Q,\partial Q).
\]
The same argument as in Lemma~\ref{lemm:H1=0} shows that the former is
isomorphic to $H^1((\cup_iM_i)/T)=H^1(\partial Q)$, and the above isomorphism
implies (through the projection $(ET\times M)/T\to M/T=Q$) that the
coboundary map $H^1(\partial Q)\to H^2(Q,\partial Q)$ in the exact sequence
of the pair $(Q,\partial Q)$ is an isomorphism.  Therefore, we get the
following exact sequence fragment:
\[
  0\to H^2(Q) \to H^2(\partial Q) \to H^3(Q,\partial Q).
\]
Since $H^2(\partial Q)\cong H^2(S^{n-1})$ by (\ref{eqn:boundary Q}), we have
$H^2(Q)=0$ if $n\ge 4$. When $n=3$, the coboundary map above is an
isomorphism because $Q$ is orientable by Lemma~\ref{lemm:H1=0}, whence
$H^2(Q)=0$ again.  This completes the proof of the claim.

\smallskip

Since $H^2(Q)=0$, we have a map $\Phi\colon M\to M_P(\L)$ as in the proof of
the \lq\lq if" part. Let us consider the diagram~\eqref{seqma} with $\k$
coefficient where $\k=\Q$ or $\Z/p$ with prime $p$. Using induction and a
Mayer--Vietoris argument, we deduce that $H^*_T(\cup_i M_{P,i};\k)\to
H^*_T(\cup_i M_{i};\k)$ is an isomorphism. We know that $H^*(P,|\P|;\k)\cong
H^*(D^n,S^{n-1};\k)$ by Lemma~\ref{kgore}, and it follows from the
construction of $\Phi$ that the induced map
\begin{equation} \label{eqn:QP}
  H^*(D^n,S^{n-1};\k)\cong H^*(P,|\P|;\k)\to H^*(Q,\partial Q;\k)
\end{equation}
is an isomorphism in degree $n$, and thus is injective in all degrees.
Therefore (an extended version of) the 5-lemma (see \cite[p.185]{Spanier})
applied to~\eqref{seqma} with $\k$ coefficient shows that $\Phi^*\colon
H^*_T(M_P;\k)\to H^*_T(M;\k)$ is injective.  Here, $H^*_T(M)\cong \Z[Q]\cong
H^*_T(M_P)$ by Corollary~\ref{coro:eqcoh} (or Proposition 3.4 in
\cite{masu99}) and Proposition~\ref{prop:srfr} (or Theorem 4.8 of
\cite{da-ja91}), so $H^*_T(M_P;\k)$ and $H^*_T(M;\k)$ have the same dimension
over $\k$ in each degree.  Therefore, the monomorphism $\Phi^*\colon
H^*_T(M_P;\k)\to H^*_T(M;\k)$ is actually an isomorphism.  Again, the 5-lemma
applied to~\eqref{seqma} with $\k$ coefficients implies that the map
(\ref{eqn:QP}) is an isomorphism, so $H^*(Q,\partial Q;\k)\cong
H^*(D^n,S^{n-1};\k)$ for any $\k$ and hence $H^*(Q,\partial Q)\cong
H^*(D^n,S^{n-1})$.  This together with (\ref{eqn:boundary Q}) (or the
Poincar\'e--Lefschetz duality) gives the acyclicity of $Q$, thus finishing the
proof of the theorem.
\end{proof}

The following statement gives a characterisation of simplicial complexes
associated with torus manifolds with cohomology generated in degree two.

\begin{theorem}
A simplicial complex $\K$ is associated with a torus manifold $M$ whose
cohomology is generated by its degree-two part if and only if $\K$ is
Gorenstein* and $\Z[\K]$ admits an l.s.o.p.
\end{theorem}
\begin{proof}
If $H^*(M)$ is generated by its degree-two part, then $\K$ is Gorenstein*, in
particular $\Z[\K]$ is Cohen--Macaulay by Lemma~\ref{kgore}.  Moreover
$H^*_T(M)\cong \Z[\K]$ by Corollary~\ref{coro:eqcoh} (or Proposition 3.4
in~\cite{masu99}). Since $H^*_T(M)\cong H^*(BT)\otimes H^*(M)$ as an
$H^*(BT)$-module by Lemma~\ref{lemm:freeness}, $\Z[\K]$ admits an l.s.o.p.

Now we prove the ``if" part. According to Theorem 12.2 of \cite{davi83},
there exists a homology polytope $Q$ whose nerve is~$\K$.  Since the face
ring $\Z[\K]$ admits an l.s.o.p., it is a free module over a polynomial ring
$\Z[t_1,\dots,t_n]$ in $n$ variables.  We can express any element $t\in
H^2(BT)\cong\Z[t_1,\dots,t_n]$ as
\[
  t=\sum_{i=1}^ma_i(t)v_i,
\]
where $a_i(t)\in\Z$. Clearly, $a_i(t)$ is linear on $t$, so $a_i$ can be
viewed as an element of the dual space $H_2(BT)$ (see
Proposition~\ref{prop:a_i}). Now define a map
$\L$~\eqref{eqn:characteristic map} by sending $Q_i$ to~$a_i$. Then
$M:=M_Q(\L)$ (see~\eqref{idspa}) is a torus manifold, and its cohomology is
generated in degree two by Theorem~\ref{hptpe}, which finishes the proof.
\end{proof}

\section{Orbit spaces of torus manifolds with vanishing odd degree cohomology}
\label{vanco}

Let $\G$ be a face of $Q$.  The facial submanifold $\MG=\pi^{-1}(\G)$ is a
connected component of an intersection of finitely many characteristic
submanifolds. The Whitney sum of their normal bundles restricted to $\MG$
gives the normal bundle $\n{\G}$ of $\MG$.  The orientations for $M$ and
characteristic submanifolds determine a $T$-invariant complex structure on
$\n{\G}$, so that the complex projective bundle $P(\n{\G})$ of $\n{\G}$ can
be considered. Replacing $\MG$ in $M$ by $P(\n{\G})$, we obtain a new torus
manifold $\til M$. The passage from $M$ to $\til M$ is called the
\emph{blowing-up} of $M$ at $\MG$. (Remark: the normal bundle $\n{\G}$ admits
many invariant complex structures and the following argument works once we
choose one.) The orbit space $\til Q$ of $\til M$ is then the result of
``cutting off" the face $\G$ from $Q$, and the simplicial cell complex dual
to $\til Q$ is obtained from that dual to $Q$ by applying a stellar
subdivision of the face dual to $\G$.

\begin{lemma}\label{qtild}
The orbit space $\til Q$ is face-acyclic if and only if so is $Q$.
\end{lemma}
\begin{proof}
By cutting the face $F$ off $Q$ we obtain a new facet $\til\G\subset\til Q$,
and all other new faces of $\til Q$ are contained in this facet. The
projection map $\til Q\to Q$ collapses $\til \G$ back to $\G$.  The face $\G$
is a deformation retract of $\til \G$.  Hence, $\G$ is acyclic if and only if
$\til \G$ is acyclic. The same is true for any other new face of $\til Q$.
It is also clear from the construction that $Q$ is a deformation retract of
$\til Q$. Therefore, $\til Q$ is acyclic if and only if so is $Q$.
\end{proof}

\begin{lemma}\label{mtild}
$H^{odd}(\til M)=0$ if $H^{odd}(M)=0$.
\end{lemma}
\begin{proof}
The facial submanifold $\MG\subset M$ is blown up to a codimension-two facial
submanifold $\til M_{\til \G}\subset\til M$, namely, $\til M_{\til \G}=
P(\n{\G})$.  Since $\til M_{\til \G}$ is the total space of a bundle with
base $\MG$ and fibre a complex projective space, its cohomology is a free
$H^*(\MG)$-module on even-dimensional generators by Dold's theorem (see,
e.g.,~\cite[Ch.~V]{ston68}).  If $H^{odd}(M)=0$, then $H^{odd}(\MG)=0$ by
Lemma~\ref{lemm:odd=0} and hence $H^{odd}(\til M_{\til \G})=0$.  Let $\til
M\to M$ be the collapse map and consider the diagram
$$
\begin{CD}
  H^{k-1}(\MG) @>>> H^k(M,\MG) @>>> H^k(M) @>>> H^k(\MG) \\
    @VVV              @VV\cong V            @VVV        @VVV     \\
  H^{k-1}(\til M_{\til \G}) @>>> H^k(\til M,\til M_{\til \G}) @>>>
   H^k(\til M) @>>> H^k(\til M_{\til \G})
\end{CD}
$$
where the second vertical arrow is an isomorphism by excision.  Assume that
$k$ is odd. If $H^{odd}(M)=0$ then $H^{k-1}(\MG)\to H^k(M,\MG)$ is onto.
Therefore, it follows from the above commutative diagram that $H^{k-1}(\til
M_{\til \G})\to H^k(\til M,\til M_{\til \G})$ is also onto.  Since $H^k(\til
M_{\til \G})=0$, this implies $H^k(\til M)=0$.
\end{proof}

The following main result of this section is an analogue of
Theorem~\ref{hptpe}.

\begin{theorem}\label{cohfa}
The odd-degree cohomology of $M$ vanishes if and only if $M$ is locally
standard and the orbit space $Q$ is face-acyclic.
\end{theorem}
\begin{proof}
The idea is to reduce to Theorem~\ref{hptpe} by blowing up sufficiently many
facial submanifolds $\MG=\pi^{-1}(\G)$.  Since the barycentric subdivision is
a sequence of stellar subdivisions, by applying sufficiently many blow-ups we
get a torus manifold $\widehat M$ with orbit space $\widehat Q$ such that the
face poset of $\widehat Q$ is the barycentric subdivision of the face poset
of $Q$. The collapse map $\widehat M\to M$ is decomposed into a sequence of
collapse maps for single blow-ups:
\begin{gather}\label{blseq}
\begin{CD}
  M=M_0 @<<< M_1 @<<< \dots @<<< M_k=\widehat M.
\end{CD}
\end{gather}

Assume that $H^{odd}(M)=0$. Then $M$ is locally standard by
Theorem~\ref{theo:local standardness}.  By applying Lemma~\ref{mtild} several
times we get $H^{odd}(\widehat M)=0$. By construction, all the intersections
of faces of $\widehat Q$ are connected, so $H^*(\widehat M)$ is generated by
its degree-two part by Theorem~\ref{theo:face CM} and $\widehat Q$ is a
homology polytope by Theorem~\ref{hptpe}. In particular, $\widehat Q$ is
face-acyclic.  Finally, by applying Lemma~\ref{qtild} inductively we conclude
that $Q$ is also face-acyclic.

The scheme of the proof of the \lq\lq if" part is same as that of
Theorem~\ref{hptpe}. But there are two things to be checked. These are
\begin{itemize}
\item[(1)] $|\P|$ has the homology of an $(n-1)$-sphere,
\item[(2)] $\Z[\P]$ is Cohen--Macaulay.
\end{itemize}
Let $\widehat \P$ be the face poset of $\widehat Q$. Since $Q$ is
face-acyclic, $\widehat Q$ is a homology polytope. Therefore, $|\widehat \P|$
has the homology of an $(n-1)$-sphere by Lemma~\ref{kgore}. However,
$|\widehat \P|=|\P|$, so the first statement above follows. Since $\widehat
Q$ is a homology polytope, $\Z[\widehat \P]$ is Cohen--Macaulay by
Lemma~\ref{kgore}. This implies that $\Z[\P]$ itself is Cohen--Macaulay by
Corollary 3.7 of \cite{stan91}, proving the second statement above.
\end{proof}

\begin{remark}
As one can easily observe, the argument in the ``only if" part of the above
theorem is independent of Theorem~\ref{theo:GKM} and
Corollary~\ref{coro:eqcoh}. Now, given that $Q$ is face-acyclic, one readily
deduces that the 1-skeleton of $Q$ is connected. Indeed, otherwise the
smallest face containing vertices from two different connected components of
the 1-skeleton would be a manifold with at least two boundary components and
thereby non-acyclic. Thus, our reference to Theorem~\ref{theo:GKM} was
actually irrelevant, although it made the arguments more straightforward.

Finally, we note that the proof of the ``if" part of Theorem~\ref{cohfa}
could have been identical to that of the ``only if" part if the converse of
Lemma~\ref{mtild} was true. It is indeed the case, however the only proof we
have so far uses quite complicated analysis of Cohen--Macaulay simplicial
posets. We are going to write it down elsewhere.
\end{remark}

\section{Gorenstein simplical posets and
Betti numbers of torus manifolds}\label{gopos}

The barycentric subdivision $\widehat \P$ of a simplicial poset $\P$ is (the
face poset of) a simplicial complex and $\P$ is called \emph{Gorenstein*} if
$\widehat \P$ is Gorenstein* (\cite{stan91}, \cite{stan96}). If $\P$ is the
simplicial poset associated with a torus manifold $M$ with $H^{odd}(M)=0$,
then the torus manifold $\widehat M$ corresponding to $\widehat \P$ has
cohomology generated by its degree-two part as remarked in the proof of
Theorem~\ref{cohfa}. Hence, $\widehat \P$ is Gorenstein* by Lemma~\ref{kgore}
and $\P$ is Gorenstein* by definition. In \cite{stan91} Stanley proved that
any vector satisfying the conditions in Theorem~\ref{theo:betti} below is an
$h$-vector of a Gorenstein* simplicial poset. He also conjectured that those
conditions are necessary. In this section we prove this conjecture for
Gorenstein* simplicial posets $\P$ associated with torus manifolds $M$ with
vanishing odd degree cohomology, and characterize $h$-vectors of those
Gorenstein* simplicial posets. The Stanley conjecture was proved in full
generality by the first author in~\cite{masu??}.

Since
\begin{equation} \label{eqn:hvector}
  h_i(\P)=\rank_\Z H^{2i}(M),
\end{equation}
by (\ref{betti}), we need to characterise the Betti numbers
of torus manifolds with vanishing odd degree cohomology. We note that
\[
  h_i(\P)\ge 0,\ h_i(\P)=h_{n-i}(\P)
  \text{\ for all $i$, and $h_0(\P)=1$.}
\]

\begin{theorem} \label{theo:betti}
Let $\mb h=(h_0,h_1,\dots,h_n)$ be a vector of non-negative integers with
$h_i=h_{n-i}$ for all $i$ and $h_0=1$. Any of the following (mutually
exclusive) conditions is sufficient for the existence of a rank $n$
Gorenstein* simplicial poset $\P$ that is associated with a $2n$-dimensional
torus manifold with vanishing odd degree cohomology and has $h$-vector~$\mb
h$:
\begin{itemize}
\item[(1)] $n$ is odd,
\item[(2)] $n$ is even and $h_{n/2}$ is even,
\item[(3)] $n$ is even, $h_{n/2}$ is odd, and $h_i>0$ for all $i$.
\end{itemize}
Moreover, if $\mb h$ is the $h$-vector of a simplicial poset of the above
described type, then it satisfies one of the above three conditions.
\end{theorem}
\begin{proof}
For a torus manifold $M$, we set $h_i(M)=\rank_\Z H^{2i}(M)$.  Thanks
to~(\ref{eqn:hvector}), we may use $h_i(M)$ instead of $h_i(\P)$ to prove the
theorem.

We shall prove the sufficiency first.  Examples~\ref{CPn} and~\ref{2nsphere}
produce torus manifolds $\C P^n$, $S^{2n}$ and $S^{2n-2k}\times S^{2k}$ with
$1\le k\le n-1$. In all three cases the odd-degree cohomology is zero.  If
$M_1$ and $M_2$ are torus manifolds (of same dimension) with vanishing odd
degree cohomology, then their equivariant connected sum $M_1\mathbin{\#}M_2$
at two fixed points with isomorphic tangential representations produces a
torus manifold with vanishing odd degree cohomology. We have
$$
  h_i(M_1\mathbin{\#}M_2)=h_i(M_1)+h_i(M_2)\qquad\text{for $1\le i\le n-1$.}
$$
Using this identity, one easily gets any vector satisfying the conditions in
the theorem by taking equivariant connected sum of $\C P^n$, $S^{2n}$ and
$S^{2n-2k}\times S^{2k}$.

Now we prove the necessity.  Let $M$ be a torus manifold of dimension $2n$.
It suffices to prove that $h_{n/2}(M)$ is even if $n$ is even and $h_i(M)=0$
for some $i>0$.

Let $G$ be the $2$-torus subgroup of $T$ of rank $n$ (that is,
$G\cong(\Z/2)^n$). Then the equivariant total Stiefel--Whitney class of $M$
with the restricted $G$-action is defined to be the ordinary total
Stiefel--Whitney class of the vector bundle $EG\times_G \mathcal TM\to
EG\times_G M$, and is denoted by $w^G(M)$.  By definition, $w^G(M)$ lies in
$H^*_G(M;\Z/2)$. We denote by $\tau_i$ the image of the identity under the
equivariant Gysin map $H^0_G(M_i;\Z/2)\to H^2_G(M;\Z/2)$, where $M_i$
$(i=1,\dots,m)$ are characteristic submanifolds of $M$.

\medskip
\noindent
{\bf Claim.} $w^G(M)=\prod_{i=1}^m(1+\tau_i)$.

\smallskip

\noindent
The proof of the claim is similar to that of Theorem 3.1 in \cite{masu99},
where the same formula was proved for the total equivariant Chern class.
Since $H^{odd}(M;\Z/2)=0$ and $M^G=M^T$ is isolated, we have
\[
  \dim H^*(M;\Z/2)=\chi(M)=\chi(M^T)=\chi(M^G)=\dim H^*(M^G;\Z/2).
\]
Therefore, $H^*_G(M;\Z/2)$ is a free $H^*(BG;\Z/2)$-module (see \cite[Theorem
VII.1.6]{bred72}).  It follows from the localisation theorem that the
restriction map
\begin{equation} \label{eqn:rest2}
  H^*_G(M;\Z/2)\to H^*_G(M^G;\Z/2)
\end{equation}
is injective. Given $p\in M^G=M^T$, set $I(p):=\{ i \colon p\in M_i\}$.  The
cardinality of $I(p)$ is $n$ and the tangential $G$-representation $\mathcal
T_pM$ decomposes as
\[
  \mathcal T_pM=\bigoplus_{i\in I(p)}\nu_i|_p
\]
where $\nu_i$ is the normal bundle of $M_i$ to $M$ and $\nu_i|_p$ is its
restriction to $p$. It follows that
\begin{equation} \label{eqn:wG}
  w^G(M)|_p=\prod_{i\in I(p)} w^G(\nu_i|_p).
\end{equation}
Since $\nu_i$ is orientable of real dimension two, $w^G_1(\nu_i)=0$ and
$w^G_2(\nu_i)$ is the mod 2 reduction of the equivariant Euler class of
$\nu_i$. Therefore, we have $w^G_2(\nu_i|_p)=\tau_i|_p$ for $i\in I(p)$.
Moreover, $\tau_i|_p=0$ for $i\notin I(p)$ by a property of equivariant Gysin
homomorphism. Thus, the identity (\ref{eqn:wG}) gives
\[
  w^G(M)|_p=\prod_{i\in I(p)}(1+\tau_i)|_p=\prod_{i=1}^m (1+\tau_i)|_p.
\]
This together with the injectivity of the restriction map in
(\ref{eqn:rest2}) proves the claim.

\medskip
The forgetful map $H^*_G(M;\Z/2)\to H^*(M;\Z/2)$ takes the equivariant
Stiefel--Whitney class $w^G(M)$ to the (ordinary) Stiefel--Whitney class
$w(M)$ of $M$.  Since $\tau_i$ is of degree two, the above claim shows that
$w_{2n}(M)$ is a polynomial in degree two elements.  Assume $h_i(M)=0$ for
some $i>0$. Then $w_{2n}(M)=0$.  The mod 2 reduction of the Euler
characteristic $\chi(M)$ of $M$ agrees with $w_{2n}(M)$ evaluated on the mod
2 fundamental class of $M$. Hence, $w_{2n}(M)=0$ implies that $\chi(M)$ is
even. Here $\chi(M)=\sum_{i=0}^n h_i(M)$ and $h_i(M)=h_{n-i}(M)$ by the
Poincar\'e duality, thus $h_{n/2}(M)$ must be even for even~$n$.
\end{proof}

\end{document}